\documentclass[12pt,a4paper]{article}
\usepackage{graphicx,bbm,color}
\usepackage{amsmath,amsfonts,amssymb,graphicx,epsfig,enumerate}
\usepackage{graphicx,bbm,color}
\usepackage{amsmath,amsthm,amsfonts,amssymb,graphicx,epsfig,enumerate}

\newcommand{\sig}{\sigma}
\newcommand{\eps}{\epsilon}
\newcommand{\rset}{\ensuremath{\mathbb{R}}}

\newcommand{\Exp}{\ensuremath{\mathbb{E}}}
\newcommand{\Hu}[2]{\mathcal{H}_{#2}\left(#1\right)}
\newcommand{\dHu}[2]{\mathcal{H}'_{#2}\left(#1\right)}

\newcommand{\IH}[2]{\mathcal{B}_{#2}\left(#1\right)}
\newcommand{\dIH}[2]{\mathcal{B}'_{#2}\left(#1\right)}
\newcommand{\ddIH}[2]{\mathcal{B}''_{#2}\left(#1\right)}
\newcommand{\indic}[1]{1\!\!1_{#1}}
\newcommand{\CA}{\mathcal{A}}
\newcommand{\CO}{\mathcal{O}}

\newcommand{\CN}{\mathcal{N}}
\newcommand{\cvp}{\stackrel{\mathbb{P}}{\rightarrow}}
\newcommand{\Pb}[1]{\mathbb{P}\left[#1\right]}

\newcommand{\X}{\ensuremath{\mathbf{X}}}
\newtheorem{theorem}{Theorem}
\newtheorem{lemme}{Lemma}

\title{The BerHu penalty and the grouped effect}

\author{Sophie Lambert-Lacroix,\\ UJF-Grenoble 1 / CNRS / UPMF / TIMC-IMAG\\ UMR 5525, Grenoble, F-38041, France
\\and Laurent Zwald,\\ LJK, Universit\'e de Grenoble et CNRS, UMR 5224\\
51, rue des Math\'ematiques, B.P. 53, 38041 Grenoble cedex 9, France}

\begin{document}

\begin{center}
{\bf The BerHu penalty and the grouped effect}
\end{center}

\begin{center}
       Sophie Lambert-Lacroix \\
UJF-Grenoble 1 / CNRS / UPMF / TIMC-IMAG\\ 
UMR 5525, Grenoble, F-38041, France\\
and\\
       Laurent Zwald\\LJK - Universit\'e de Grenoble\\
BP 53, 38041 Grenoble cedex
9, France
\end{center}

{\bf Abstract.}  

The Huber's criterion is a useful method for
robust regression. The adaptive least absolute shrinkage and
selection operator (lasso) is a popular technique for simultaneous
estimation and variable selection. In the case of small sample size and large covariables numbers, this penalty is not very satisfactory variable selection method. In this paper, we introduce an adaptive reversed version of Huber's criterion as a penalty function. We call this penalty adaptive Berhu penalty. As for elastic net penalty, small coefficients contribute their $\ell_1$ norm to this penalty while larger coefficients cause it to grow quadratically (as ridge regression). We show that 
the estimator associated with criterion such that ordinary least square or Huber's one combining with adaptive Berhu penalty enjoys the oracle properties. 
In addition, this procedure encourages a grouping effect.
This approach is compared with adaptive elastic net regularization. Extensive simulation studies demonstrate satisfactory finite-sample performance of such procedure.
A real example is analyzed for illustration purposes.


{\bf Availability.} The software that implements the
procedures on which this paper focuses is developed in {\tt Matlab}. It is available at {\tt
http://ljk.imag.fr/membres/Laurent.Zwald}.

\section{Introduction}

Data subject to heavy-tailed errors or outliers are commonly encountered in applications which may appear either in response variables or in the predictors. We consider here the regression problem with eventually responses subject to heavy-tailed errors or outliers. In this case, the Ordinary Least Square (OLS)  estimator  is reputed to be not efficient. To overcome this problem, the least absolute deviation (LAD) or Huber type estimator for instance can be useful.
On the other hand, an important topic in linear regression analysis is variable selection. Variable selection is particularly important when the true underlying model has sparse representation. To enhance the prediction performance of the fitted model and get an easy interpretation of the model, we need to identify significant predictors. Scientists prefer a simpler model because it puts more light on the relationship between the response and covariates. We consider the important problem of robust model selection.

The lasso penalty is a regularization technique for simultaneous estimation and variable selection (\cite{Tib96}). It consists to introduce $\ell_1$ penalty. This penalty forces to shrink some coefficients. In \cite{Fan01}, the authors show that since lasso uses the same tuning parameters for all the regression coefficients, the resulting estimators may suffer an appreciable bias. Moreover in the case of the small sample $n$ with larger number of covariables $p$, the lasso selects at most $n$ variables. 
Recently, \cite{MeinshausenBuhlman04,LengLinWahba06,ZhaoYu06}  and \cite{Zou06} show that the underlying model must satisfy a nontrivial condition for the lasso estimator be consistent in variable selection. Consequently, in some cases, lasso estimator cannot be consistent in variable selection.  For instance,  \cite{Zou06}  assigns adaptive weights for penalizing differently coefficients in the $\ell_1$ penalty and calls this new penalty the adaptive lasso. These adaptive weights in the penalty allow to have the oracle properties. Moreover, the adaptive lasso can be solved by the same efficient algorithm (LARS) for solving lasso (see \cite{Zou06}). Notice that recently (see \cite{Lam10}), this penalty has been combined with Huber's criterion. The estimator associated with this procedure enjoys oracle properties. 

On the other hand, if there is a group of variables among which the pairwise correlations are very high, then the lasso penalty tends to select only any one variable from this group. Ridge regression ($\ell_2$ penalty) does not make variables selection but tends instead to share the coefficients value among the group of correlated predictors. Moreover if there exist high correlations among predictors, the prediction performance of ridge regression dominated the lasso \cite{Tib96}.
In order to overcome to this drawback of the lasso, \cite{Zou05} proposes a new regularization technique that combines the lasso and the ridge penalties. They call their method ``elastic net'' (en). The en penalty is the sum of the lasso and the ridge penalties. 
However even for usual case, it does not deemed to be an oracle procedure. In \cite{Gosh2007}, the author proposes a new version of the elastic net called adaptive elastic net (adaptive en) which inherits some of the desirable properties of the adaptive lasso and elastic net. He proves its oracle properties.
In \cite{Owe06}, the author proposes to use a reversed version of Huber's criterion (called Berhu) as a penalty function. Let us recall that the Huber criterion (see \cite{Hub81}) is a hybrid of squared error for relatively small errors and absolute error for relative large ones. The Berhu penalty is such that relatively small coefficients contribute their $\ell_1$ norm to this penalty whiles larger ones cause it to grow quadratically. This hybrid sets some coefficients to 0 as the lasso does while shrinking the larger coefficients in the same way as ridge regression. In \cite{Owe06}, the author provides some way in order to optimize some objective function constituted of both the Huber criterion and the Berhu penalty in a no-adaptive form. Nevertheless nothing is shown about asymptotic feature.  

In this paper we introduce an adaptive Berhu penalty with concomitant. We use it with the ordinary least square criterion or the Huber's one in order to take into account of data subject to heavy-tailed errors or outliers. We show that the estimator associated with such procedures enjoys the oracle properties (in the standard case of least square criterion and in the case of the Huber's one).  In addition this procedure encourages a grouping effect in the following way.  The spirit of the Berhu penalty with concomitant implicitly is to create one group with the largest coefficients. This group is penalized in a $\ell_2$ way like the grouped lasso of \cite{YuanLin2006} to avoid to remove anyone of these largest coefficients. The smallest coefficients are treated individually by an $\ell_1$-penalty. 
The en procedure relies on the fact that, in order to have a grouped effect, we want to keep or delete together high correlated variables. We show that when combining with ordinary least squares criterion, the Berhu penalty leads to this ``grouping effect property''.

The rest of the article is organized as follows. In Section~2, we introduce the adaptive BerHu penalty and show that it induces a grouped effect. In Section~3, we give its statistical properties. Section~4 is devoted to simulation and illustration over real data. This study compares the least square criterion and the Huber's criterion with  various penalties such as adaptive lasso, ridge, en and adaptive Berhu.  All  technical proofs are relegated to the Appendix.

\section{The Berhu penalty}

\subsection{The adaptive Berhu}\label{sec:adaber}

Let us consider the linear regression model
\begin{equation}\label{eq:model}
y_i = \alpha^* + \mathbf{x}_i^T\beta^* + \sigma\epsilon_i,\quad i=1,\ldots, n,
\end{equation}
where $\mathbf{x}_i = (x_{i1},\ldots,x_{ip})^T$ is the
$p$-dimensional centered covariable (that is $\sum_{i=1}^n \mathbf{x}_i=0$), 
$\alpha^*$ is the constant parameter and
$\beta^*=(\beta^*_1,\ldots,\beta^*_p)^T$ are the associated regression
coefficients. We suppose that $\sigma>0$ and  $\epsilon_i$ are
independent and identically-distributed random errors with mean 0
and variance 1, when it exists. Indeed in the sequel we do not need existence of variance.
Let $\CA=\{1\leq j\leq p,\,\beta^*_j\neq0\}$ and $p_0=\vert\CA\vert$. In variables selection context, we usually assume
that $\beta^*_j\neq 0,$ for $j\leq p_0$ and  $\beta^*_j= 0,$ for
$j>p_0$ for some $p_0\geq 0.$ In this case the correct model has
$p_0$ significant regression variables. We denote by $\beta_{\CA}$ the vector given by the coordinates of $\beta$ the index of which are in $\CA$. 

When $p_0=p,$ the unknown parameters in the model (\ref{eq:model})
are usually estimated by minimizing the ordinary least squares criterion. To shrink unnecessary coefficients to 0,
\cite{Tib96} proposed to introduce a constraint on the $\ell_1$-norm of the coefficients: 
$$\sum_{i=1}^n(y_i-\alpha-\mathbf{x}_i^T\beta)^2+\lambda_n\sum_{j=1}^p|\beta_j|\,,$$
where $\lambda_n>0$ is the tuning parameter. Notice that the intercept $\alpha$
does not appear in the penalty term since it is not reasonable to constrain it.

Lots of reproaches have already been done to the Lasso  (see e.g. \cite{Zou05}). In this paper, we focuse on the fact that when some variables are highly correlated, the $\ell_1$ penalty  tends to keep only one variable for each group. The literature already contains attempts to solve this problem. To begin with, grouped lasso procedures have been proposed first in \cite{YuanLin2006} where the $\ell_1$ penalty is imposed on predefined groups of coefficients. More precisely, the penalty is the $l_1$-norm of the vector composed of the $\ell_2$-norm of each group of coefficients:
$$
\sum_{i=1}^n(y_i-\alpha-\mathbf{x}_i^T\beta)^2+\lambda_n\sum_{j=1}^L\sqrt{p_l}\|(\beta)_j\|_{2},
$$
where $(\beta)_j$ is the coordinates bloc corresponding to the $j$-th group. Consequently, the sparsity is encouraged at the group level (see also \cite{ZhaopengRochaGuilhermeYuBin2009} and \cite{HastieTibshiraniFriedman2003} page 91 for further references). In our framework it is diffcult to use the approach of  group lasso since there is no obvious way for choosing the groups a priori. 
Next, \cite{Zou05} has proposed the Elastic Net. The naive Elastic Net is obtained by minimizing:
\begin{equation}
\label{EN}
\sum_{i=1}^n(y_i-\alpha-\mathbf{x}_i^T\beta)^2+\lambda_{1,n}\sum_{j=1}^p|\beta_j|+\lambda_{2,n}\sum_{j=1}^p\beta_j^2,
\end{equation}
and the Elastic Net is a modification of this.
In this procedure, the penalty imposed on the small coefficients is the sum of an $\ell_1$-norm and a squared $\ell_2$-norm. Moreover, ridge penalty reduces the variance of the estimates by imposing a small squared norm of all the coefficients. However, it suffices to constraint the largest coefficients  to be small to get this reduction of variance: by definition, the smallest one do not need to be constrained to be small. Consequently, we consider a penalty which is quadratic only on the largest coefficient. Following \cite{Owe06}, we focused on a penalty that acts separately on small and large coefficients. We consider the Berhu penalty defined by 
\begin{equation}
\label{defBerHu}{\cal B}_L(z)=\left\{
\begin{array}{ll}
|z|&|z|\leq L,\\
\frac {z^2+L^2}{2L}&|z|> L,
\end{array}
\right.\end{equation}
where $L$ is any positive real. As Huber criterion, the Berhu function needs to be scaled. Precisely, the penalty can be defined by
$$
\sum_{j=1}^p {\cal B}_L\left(\frac {\beta_j} \tau \right),
$$
where $\tau$ is a scale parameter to be determined.  To  do that we can as in  \cite{Owe06} replace the penalty term by
$$\text{pen}(\beta)=\min_{\tau >0}\left(p\tau+\tau\sum_{j=1}^p{\cal B}_L\left(\frac{\beta_j}{\tau}\right)\right).$$

Fan and Li \cite{Fan01}  showed that the lasso method leads
to estimators that may suffer an appreciable bias. Furthermore
they conjectured that the oracle properties do not hold for the
lasso. Hence Zou \cite{Zou06} proposes to consider the following
modified lasso criterion, called adaptive lasso,
$$
\sum_{i=1}^n(y_i-\alpha-\mathbf{x}_i^T\beta)^2+\lambda_n\sum_{j=1}^p \hat{w}_j^{adl} |\beta_j|,$$
where $ \mathbf{\hat{w}}^{adl}=(\hat{w}_1^{adl},\ldots,\hat{w}_p^{adl})$ is a known weights vector.
This modification allows to produce sparse solutions more effectively than lasso. Precisely, Zou \cite{Zou06}
shows that with a proper choice of $\lambda_n$ and of $\mathbf{\hat{w}}^{adl}$ 
the adaptive lasso enjoys the oracle properties. Such a penalty has been used in the en penalty (see \cite{Gosh2007}).

Here we propose to make the Berhu penalty adaptive. That is we consider the following penalty $\min_{\tau\in \rset} P^{adb}(\beta,\tau)$ with
$$P^{adb}(\beta,\tau)=\left\{\begin{array}{lll}
\tau \left(\sum_{j=1}^{p} \frac 1 {\hat{w}_j^{adb}}+ \sum_{j=1}^p \hat{w}_j^{adb}{\cal B}_L\left(\frac {\beta_j} \tau\right)\right)&if&\tau > 0,\\
0& if&  \beta=0,\; \tau=0,\\
+\infty& if&\beta\neq 0,\;  \tau=0.
\end{array}
\right.
$$
where $ \mathbf{\hat{w}}^{adb}=(\hat{w}_1^{adb},\ldots,\hat{w}_p^{adb})$ is a known weights vector.
We will see at Section 3 that the resulting estimator enjoys the oracle properties. Let us notice that \cite{Owe06} introduced the Berhu penalty in his no-adaptive form and in the context of robust regression only. Moreover nothing is shown about asymptotic feature.

In the general case, the (adaptive) Behru penalty behaves like lasso on the smallest coefficients and does not delete the largest ones, whatever the correlation structure. That can be what we except to a right model selection procedure. This interpretation relies on the following calculation when $\beta$ is fixed. For instance in the non adaptive case, let us sort the absolute values of the coordinates of $\beta$:
$$
\vert\beta_{(p)}\vert\leq\cdots\leq\vert\beta_{(1)}\vert\,.
$$
Let $k(\beta)$ denote the number of non-zeros coefficients of $\beta$. Then the minimum defined in $\text{pen}(\beta)$ is achieved at 
$$\hat{\tau}(\beta)=\sqrt{\frac{1}{2Lp+L^2(q(\beta)-1)}\sum_{j=1}^{q(\beta)-1}\beta_{(j)}^2},$$
if $\beta\neq0$ and where $q(\beta)$ is the unique integer between 2 and $k(\beta)+1$ such that ${\vert\beta_{(q(\beta))}\vert}/{L}\leq\hat{\tau}(\beta)\leq{\vert\beta_{(q(\beta)-1)}\vert}/{L}$. Consequently,
\begin{equation}
\label{exprpen}
\text{pen}(\beta)=\sqrt{\frac{2p}{L}+q(\beta)-1}\sqrt{\sum_{j=1}^{q(\beta)-1}\beta_{(j)}^2} + \sum_{j=q(\beta)}^{k(\beta)}\vert\beta_{(j)}\vert\ .
\end{equation}
The en procedure (or its variant Elastic Corr-Net \cite{Anbari2008}) relies (explicitly for Elastic Corr-Net) on the fact that, in order to have a grouped effect, we want to keep or delete together high correlated variables. We will see that it is the case for Berhu procedure in Section~\ref{sec:groupinfeffect}. But we can note here different spirit of the Berhu penalty with concomitant: it implicitly creates one group with the largest coefficients (see (\ref{exprpen})). This group is penalized in a $\ell_2$ way like the grouped lasso of \cite{YuanLin2006} to avoid to remove anyone of these largest coefficients. Let us note that as in the grouped lasso penalty, the $\ell_2$-norm of the $q(\beta)-1$ largest coefficients is scaled by the squared root of the number of such coefficients present in this group. The smallest coefficients are treated individually by an $\ell_1$-penalty (see (\ref{exprpen})). Consequently, whatever the structure of the correlation matrix, the Berhu penalty with concomitant tends to keep all the largest coefficients and to delete the smallest ones. 

\subsection{Robust estimation}

To be robust to the heavy-tailed errors or outliers in the
response, a possibility is to use the Huber's criterion as
loss function as introduced in \cite{Hub81}. For any positive real $M,$ let us
introduce the following function
$${\cal H}_M(z)=\left\{
\begin{array}{ll}
z^2&|z|\leq M,\\
2M|z|-M^2&|z|> M.
\end{array}
\right.$$ This function is quadratic in small values of $z$ but
grows linearly for large values of $z$. The parameter $M$
describes where the transition from quadratic to linear takes
place. The Huber's Criterion with concomitant scale defined by
$${\cal L}_{{\cal H}}(\alpha,\beta,s)=\left\{
\begin{array}{lll}
ns + \sum_{i=1}^n {\cal H}_M\left(\frac{y_i-\alpha-\mathbf{x}_i^T\beta}s\right)s&\text{if}&s>0,\\
2M\sum_{i=1}^n\vert y_i-\alpha-\mathbf{x}_i^T\beta\vert&\text{if}& s=0,\\
+\infty &\text{if}& s<0,
\end{array}
\right.
$$
which are to minimize with respect to $s\geq0$, $\alpha$ and $\beta.$ 
To get a robust scale invariant Lasso type procedure, \cite{Lam10} proposes to minimize simultaneously over $s,\alpha$ and $\beta$ the function 
\begin{equation} 
\label{Oppaper1}
{\cal L}_{{\cal H}}(\alpha,\beta,s)+\lambda_n\sum_{j=1}^p \hat{w}_j^{adh}|\beta_j|.
\end{equation}
where $ \mathbf{\hat{w}}^{adh}=(\hat{w}_1^{adh},\ldots,\hat{w}_p^{adh})$ is a known weights vector. The loss function involving a concomitant estimation of the scale and location parameter was first proposed by Huber (\cite{Hub81}). We propose here to use the concomitant estimation of Huber with the Berhu penalty:
\begin{equation}
\label{OpIH}
Q^{{\cal H}adb}(\alpha,\beta,s,\tau)={\cal L}_{{\cal H}}(\alpha,\beta,s)+\lambda_n 
P^{adb}(\beta,\tau).
\end{equation}
This criterion is minimized simultaneously over $\alpha\in \rset,\beta\in\rset^p,s\in\rset_+$ and $\tau\in\rset_+$. So we get another scale invariant robust location estimation. Contrary to the procedure proposed in \cite{Lam10}, the largest coordinates of $\beta$ are quadratically penalized.  

\subsection{Tuning parameter estimation}\label{ss:hyperparameter}

Let us now consider the problem of tuning parameter estimation. To run these procedures we have to determine the weights vector in the adaptive penalties, the regularization constant $\lambda_n$, the parameter $M$ for Huber's criterion and  $L$ for Berhu's penalty. Usually the weights vector is given by (see  \cite{Zou06,Lam10})  $\hat{{w}}_j^{adl}=|\hat{\beta}^{unpen}_j|^{-\gamma},$ $j=1,\ldots,p$, where $\gamma>0$ and $\hat{\beta}^{unpen}$ denotes the unpenalized estimator. For instance, in the least squares context $\hat{\beta}^{unpen}$ is the ordinary least squares estimator. In fact this estimator only must be root-$n$-consistent estimator of $\beta^*$. Let us note that the theoretical part is given for these forms of  weights vector and $\gamma$ is fixed to be equal to 1
for the numerical results. For Huber's Criterion with concomitant scale we need value for $M$. As in \cite{Hub81}, we fix $M=1.345$. 
For Berhu's penalty we fix as in \cite{Owe06}, $L=M.$ Let us note that we do not have any justification to do that. However in practice we have observed that these parameters have little impact on the results.

To find optimal values for $\lambda_n,$ we use BIC-type criterions. When using least squares criterion we consider the classical BIC criterions (\cite{Schwarz1978}), That is it is recommended to select $\lambda_n$ minimizing  
$$\log\left(\sum_{i=1}^n\left(y_i-\widehat{\alpha}_{\lambda_n}-\mathbf{x}_i^T\widehat{\beta}_{\lambda_n}\right)^2\right)+k_{\lambda_n}\frac{\log(n)}{n},$$ 
over $\lambda_n,$  where $k_{\lambda_n}$ denotes the model dimension. Following  \cite{WangLeng2007} and \cite{WangLiTsai2007}, we determine $k_{\lambda_n}$ by the number of non-zero coefficients of the estimator.  When using Huber's criterion, we consider the BIC-type procedure introduced in \cite{Lam10}: we select $\lambda_n$ by minimizing 
$$\log\left({\cal L}_{{\cal H}}\left(\widehat{\alpha}_{\lambda_n},\widehat{\beta}_{\lambda_n},\widehat{s}_{\lambda_n}\right) \right)+ k_{\lambda_n}\frac{\log(n)}{2n},$$ 
over $\lambda_n$. As previously,  $k_{\lambda_n}$ denotes the number of non-zero coefficients of $\widehat{\beta}_{\lambda_n}.$ 

\subsection{The Berhu penalty with concomitant induced a grouped effect}\label{sec:groupinfeffect}

An algorithm is said to satisfy the grouping effect property if high correlated variables lead to similar estimations of the corresponding coefficients. Such a property was a motivation to introduce the Ridge Regression (\cite{HoerlKennard1970}). Indeed, the 
normal equations associated to the Ordinary Least Square do not imply any stability of the coefficients associated to highly correlated variables. Now, adding a  squared $\ell_2$-norm penalty, the corresponding normal equations imply a stability 
of the coefficients associated to highly correlated variables. Such a reasoning leads to a bound quantifying the grouping effect of the Elastic Net (\cite{Zou05}).  Such a property was generalized 
to the adaptive Elastic Net in \cite{Gosh2007} and also proved for the algorithm of \cite{BondellreichOscar2008}.

The goal of the following theorem is to provide a quantitative description for the grouping effect of the Berhu penalty with concomitant. 
\begin{theorem}\label{theo:groupingeffect}
Let $\gamma>0$ and  ($\hat{\alpha}^{adb},\hat{\beta}^{adb},\hat{\tau}^{adb}$) be a minimizer of
$$
\sum_{i=1}^n\left(y_i-\alpha-\mathbf{x}_i^T\beta\right)^2+\lambda_n P^{adb}(\beta,\tau),$$
over $\alpha\in \rset,\beta\in\rset^p$ and $\tau\in\rset_+$. We suppose that $\lambda_n>0$, $\hat{\beta}^{adb}_i\neq0$, $\hat{\beta}^{adb}_j\neq0$. In this situation, 
the following bound holds:
\begin{equation}
\label{to}
\vert
\hat{\beta}_i^{adb} \hat{w}_i^{adb}-\hat{\beta}_j^{adb} \hat{w}_j^{adb}\vert
\leq \frac{2L\hat{\tau}}{\lambda_n}\vert\vert \underline{y} \vert\vert_2\sqrt{\vert\vert x_i\vert\vert_2^2+\vert\vert x_j\vert\vert_2^2-2C_{i,j}{x_{i}}^Tx_{j}}
\end{equation}
where
$C_{i,j}=\min\left(1,\frac{\vert\hat{\beta}^{adb}_j\vert}{L\hat{\tau}^{adb}},\frac{\vert\hat{\beta}^{adb}_i\vert}{L\hat{\tau}^{adb}},\frac{\vert\hat{\beta}^{adb}_i\hat{\beta}^{adb}_j\vert}{(L\hat{\tau}^{adb})^2}\right)$.
\end{theorem}
To obtain this result for Huber's loss is a difficult task. That is an open question that is left for future work. Let us remark that when the variables are standardized in $\ell_2$-norm, this leads to 
\begin{equation}
\label{gref}
\vert
\hat{\beta}^{adb}_i \hat{w}_i^{adb}-\hat{\beta}^{adb}_j \hat{w}_j^{adb}\vert
\leq \frac{2L\hat{\tau}}{\lambda}\vert\vert \underline{y} \vert\vert_2\sqrt{2\left(1-C_{i,j}{x_{i}}^Tx_{j}\right)}\,.
\end{equation}
With $\gamma=0$, we exactly get the grouping effect property in the non-adaptive case. Let now $\gamma\in \rset_+$. The upper bound of equation \eqref{gref} is a decreasing function 
of the correlation ${x_{i}}^Tx_{j}$ between variables $i$ and $j$ (since $C_{i,j}>0$). To ensure that the coefficients $\hat{\beta}^{adb}_i$ and $\hat{\beta}^{adb}_j$ become similar if the correlation increases, from (\ref{to}), the initial 
estimator $\hat{\beta}^{unpen}$ used in the weights  $\hat{w}_i^{adb}$ has to satisfy the grouping effect property. Consequently, this bound effectively provides a quantitative description for the grouping effect of the Berhu penalty with concomitant if, for example, 
the initial estimator is obtained with a ridge penalty.

As compared with the Elastic Net bounds provided by \cite{Zou05} and \cite{Gosh2007}, we do not have to suppose that $\hat{\beta}^{adb}_i$ and $\hat{\beta}^{adb}_j$ have the same sign. Moreover, in our case, the grouping effect occurs more accurately for 
large coefficients (see Section~\ref{sec:adaber}) which is the natural situation where it has to happen. For the adaptive elastic net,  \cite{Gosh2007} also have to suppose that the initial estimator satisfies the grouping effect property.
 Moreover, \cite{ZouZang2009} recomands to choose a non-adaptive elastic net estimator as an initial estimator in the weights of the adaptive elastic net.

In the simulation study below, the initial estimator used for the weights of all adaptive methods is the corresponding unpenalized estimator. This choice avoids choosing a supplementary parameter 
(e.g. the regularization parameter of ridge regression) and also avoids numerical problems du to too small coefficients of the initial estimator. This unpenalized parameter does not satisfy the grouping effect property but comparisons between various methods remains fair. 
Moreover, in the simulation studies involving the Berhu penalty with concomitant, the variables were not normalized in $\ell_2$-norm. Indeed, using the way we get the design matrix $\X$, explicit calculations when variables  
are normalized or not leads to the same order for the corresponding upper bounds.

\section{Oracle Properties}
\label{s:theorie}

In this section we give the asymptotic properties of 
the concomitant estimator of Huber with the Berhu penalty. We show that it
enjoys the oracle properties. We have the same property by replacing  Huber's loss
by least squares one's. When necessary, we give the difference (for example for the assumptions) 
between the two loss functions.

Let $\X$ denotes the design matrix i.e. the $n\times p$ matrix the
$i^{th}$ rows of which is $\mathbf{x}_i^T$. We will use some of the
following assumptions on this design matrix.
\begin{description}
\item[(D1)]$\max_{1\leq i \leq n} {\|\mathbf{x}_i\|}/{\sqrt{n}}
\rightarrow0 \,\,\text{as} \,\, n\rightarrow\infty$\,.
\item[(D2)]${\X^T\X}/{n}\rightarrow V \,\,\text{as} \,\, n\rightarrow
\infty$ with $V_{1,1}>0,$ where $V_{1,1}$ is the first $p_0\times
p_0$ bloc of $V$, corresponding to the covariables associated with non
zero coefficients.
\end{description}
Assumption {\bf(D1)} and {\bf(D2)} are classical.  It can 
be seen as a ``compacity assumption'': it is satisfied if the variables
are supposed to be bounded. When considering  least squares criterion as loss function, we 
need only the assumption {\bf(D2)} (see for example \cite{Zou06}) while considering Huber's criterion we need the both {\bf(D1)} and {\bf(D2)} (see \cite{Lam10}).

Let us denote by $\eps$ a variable with the same law as $\eps_i$,
$i=1,\ldots, n.$ As in \cite{Lam10}, we define
$$
s^*=\underset{s>0}{\text{argmin}}\,F(s),
$$
where for $s>0$,
$$
F(s)=\Exp\left[\frac{1}{n}{\cal L}_{{\cal H}}(\alpha^*,\beta^*,s)\right
]=s+s\Exp\left[{\cal H}_M\left(\frac{\sigma\epsilon}{s}\right)\right].
$$
In addition, let us define $\tau^*>0$ satisfying 
\begin{equation}
\label{deftaustar}
\tau^*=\underset{\tau>0}{\text{argmin}}\,
\tau \left( \sum_{j=1}^{p}\vert\beta^*_j\vert^{\gamma}+\sum_{j=1}^{p_0}\frac{1}{\vert\beta^*_j\vert^{\gamma}}\IH{\frac{\beta^*_j}{\tau}}{L}\right).
\end{equation}
The following assumptions on the errors are used in
the following:
\begin{description}
\item[(N0)] The distribution of the errors does not charge the points
$\pm Ms^*$:
\[
\Pb{\sig\eps=\pm Ms^*}=0.
\]
\item[(N1)] The variable $\epsilon$ is symmetric (i.e. $\epsilon$ has
the same distribution as $-\epsilon$).
\item[(N2)] For all $a>0$, $\mathbb{P}\left[\epsilon\in[-a,a]\right]>0\,
. $
\end{description}
Note that {\bf(N0)} holds if $\eps$ is absolutely continuous with
respect to the Lebesgue's measure and {\bf(N2)} is satisfied if,
moreover, the density is continuous and strictly positive at the origin
(which is assumption A of \cite{Wan07}). Condition {\bf(N1)} is
natural without prior knowledge on the distribution of the errors and
{\bf(N2)} ensures that the noise is not degenerated. It is noticeable
that there is no integrability condition assumed on the errors $\eps$.
These three assumptions stand for the Huber's loss. For the penalized least squared
estimators (e.g. \cite{KnightFu2000} and \cite{Zou06}) we assume
that $\eps_i$ are independent identically distributed random variables 
with mean 0 and has a finite variance.

Let $(\hat{\alpha}^{adb},\hat{\beta}^{adb},\hat
{s}^{{\cal H}adb},\hat{\tau}^{{\cal H}adb})$ be defined by the minimizer of $Q^{{\cal H}adb}(\cdot)$ where 
$\hat{w}^{adb}_j$ $={1}/{\vert\hat{\beta
}^{unpen}_j\vert^{\gamma}}$ with $\gamma>1/3$ and $\hat{\beta}^{unpen}$ a
root-$n$-consistent estimator of $\beta^*$ (i.e. $\sqrt{n}(\hat{\beta
}-\beta^*)=\CO_P(1))$. We denote $\CA_n=\{1\leq j\leq p,\, \hat{\beta
}^{{\cal H}adb}_{j}\neq0\}$. Let us remark that if $\lambda_n>0$, the
argminimum $(\hat{\alpha}^{{\cal H}adb},\hat{\beta}^{{\cal H}adb},\hat
{s}^{{\cal H}adb},\hat{\tau}^{{\cal H}adb})$ exists since the criterion $Q^{{\cal H}adb}(\cdot
)$ is a convex and coercive function. 

In the following theorem we show that, with a proper choice of $\lambda
_n$, the proposed estimator enjoys the oracle properties. Its proof is
postponed in Appendix~\ref{approofthe2}.
\begin{theorem}
\label{8} Suppose that $\lambda_n/n^{\gamma\wedge 1/2}\rightarrow0,$  $\lambda
_n n^{(\gamma-1)/2}\rightarrow\infty,$ $\lambda
_n \rightarrow\infty$ and $\lambda_n>1/3$. Let us also assume that
conditions $M>1$, $p_0>0$, {\bf(N0)}, {\bf(N1)}, {\bf(N2)}, {\bf
(D1)} and {\bf(D2)} hold. Moreover, for $j=1,\ldots,p$, the weights in
$Q^{{\cal H}adb}$ are $\hat{w}^{adb}_j={1}/{\vert\hat{\beta
}^{unpen}_j\vert^{\gamma}}$ where $\hat{\beta}^{unpen}$ is a root-n-consistent
estimator of $\beta^*$. Then, any minimizer $(\hat{\alpha}^{{\cal H}adb},\hat{\beta}^{{\cal H}adb},\hat{s}^{{\cal H}adb},\hat{\tau}^{{\cal H}adb})$ of $Q^{{\cal H}adb}$ satisfies the following:
\begin{itemize}
\item Consistency in variable selection: $\mathbb{P}\left[\CA_n=\CA
\right]\rightarrow1$ as $n\rightarrow+\infty$.
\item Asymptotic normality:
\[
\sqrt{n}\left(\hat{\alpha}^{{\cal H}adb}-\alpha^*,\hat{\beta}_{\CA
}^{{\cal H}adb}-\beta^*_{\CA},\hat{s}^{{\cal H}adb}-s^*, \frac {\sqrt{\lambda_n}}{\sqrt{n}}(\hat{\tau}^{{\cal H}adb}-\tau^*)\right)
\longrightarrow_d\CN_{p_0+3}\left(0,\Sigma^2\right),
\]
where $\Sigma^2$ is the squared block diagonal matrix
\[
\Sigma^2={\tt diag}\left(\frac{\Exp\left[\dHu{\frac{\sig\epsilon
}{s^*}}{M}^2\right]}{4A_{s^*}^2},\frac{\Exp\left[\dHu{\frac{\sig\epsilon
}{s^*}}{M}^2\right]}{4A_{s^*}^2}V_{1,1}^{-1},\frac{\Exp\left[Z^2\right
]}{4D_{s^*}^2},0\right)\,
\]
and where
\[
D_{s^*}=\frac{1}{s^{*3}}
\Exp\left[\sig^2\eps^2\indic{\vert\sig\eps\vert\leq Ms^*}\right],\;
A_{s^*}=\frac{1}{s^*}\Pb{\vert\sig\eps\vert\leq Ms^*},
\]
\[
Z=1+\Hu{\frac{\sig\eps}{s^*}}{M}-\frac{\sig\eps}{s^*}\dHu{\frac{\sig\eps
}{s^*}}{M}.
\]
\end{itemize}
\end{theorem}
Analogous results hold for the least squares loss function. In this case
($M=+\infty$), the asymptotic variance
matrix ${\Exp[\dHu{\sig\eps}{Ms}^2]}V_{1,1}^{-1} /(4A_{s^*}^2)$
obtained in Theorem~\ref{8} is equal to $\sigma^2V_{1,1}^{-1} $ and we
find the asymptotic variance of theorem 2 of~\cite{Zou06}.

\section{Some numerical experiments}
In this section, we consider the both criterions least squares and Huber's one combined with the following penalties: adaptive lasso, ridge, adaptive en and adaptive Berhu. We call these methods respectively {\tt ad-lasso}, {\tt ridge}, {\tt ad-en},  {\tt ad-Berhu} , {\tt Huber-ad-lasso},  {\tt Huber-ridge},  {\tt Huber-ad-en} and {\tt Huber-ad-Berhu}. The adaptive weights are obtained from the corresponding unpenalized estimator and $\gamma=1$. 

\subsection{Simulation Results}
Here our aim is to compare the finite sample performances of these procedures. 
Paragraph \ref{mod} presents the studied models. The way simulations are conducted is described in \ref{pred} and an insight of conclusions is provided in paragraph \ref{comp}.

\subsubsection{Models used for simulations}
\label{mod}
The models used to compare the performances of the algorithms are inspired by those presented in \cite{Zou05}. They involve groups of highly correlated variables: the block-variables model (\cite{Zou05}, {\it example 4}). Let us remark that \cite{Zou05} considered a model without intercept. We now recall the definition of this model in a different way. Our formulation allows to clearly identify the groups of influencing correlated variables. They all have the form 
$
\underline{y}=\indic{n}+\X\beta^*+\sig\underline{\eps}\,,
$
where $\indic{n}$ denotes the vector of $\rset^n$ composed of ones and $\underline{y}$ (resp. $\underline{\eps}$) represents the response (resp. error) vector $(y_1,...,y_n)^T$ (resp. $(\eps_1,...,\eps_n)^T$). 
The design matrix $\X$ is constructed as follows. The rows of  $\X$ are given by $n$ independent gaussian vectors  ${\cal N}_{40}(0,\Sigma)$. They are  normalized such that the corresponding $p$-dimensional covariables are centered (as assumed in \eqref{eq:model}). The variance matrix of the variables is a block diagonal matrix of size $40$. The first block is the squared matrix of size $5$ composed of $1$ outside the diagonal and taking values $1.01$ on the diagonal. The second and third blocks are the same as the first one. The last block is the identity matrix of size $25$. The vector of true coefficients $\beta^*$ is defined as follows: the $15$ first coordinates are equal to $3$ and the $25$ last coefficients are $0$. This means that, in this model, only the $15$ first variables are influencing the response. The $25$ others are pure noise. Amongst the $15$ influencing variables, there is three groups of highly correlated variables: these groups are composed of the first five variables, the next five ones and the five last ones. The variables of different groups are independent. As compared with \eqref{eq:model}, this means that the intercept of the model is $\alpha^*=1$ and the number of variables (without the intercept) is $p=40$. Depending on the nature of the noise, various models are considered. 

\begin{itemize}
\item Model 1: {\it block-variables model, gaussian noise}. In this case,  the standard deviation of the noise is $\sig=15$ and the variables $\eps_1,\cdots,\eps_n$ are independent standard normal variables. Except for the part of the intercept parameter, this exactly example 4 of \cite{Zou05}.
\item Model 2: {\it block-variables model, mixture of gaussians}. In this case, the variables $\eps_1,\cdots,\eps_n$ are independent mixture of gaussians. Precisely, with probability $0.9$, $\eps$ is a standard normal variable and with probability $0.1$, $\eps$ is a centered normal with variance $225$.  The value $\sig=3.1009$ has been chosen such that the standard deviation of the noise is the same as in model 5. The common value is $\text{std}(\sig\eps)=3.1009\sqrt{1+0.1(225-1)}=15$.
\item Model 3: {\it block-variables model, double-exponential noise}. $\eps=D/\sqrt{\text
{var}(D)}$ and $\sigma=10.6$. The distribution of $D$ is a
standard double exponential i.e. its density is $x\in\rset\rightarrow
e^{-\vert x \vert}/2$ and $\text{var}(D)=2$.
\end{itemize}

These three models create a grouped variables situation. They allow us to illustrate the grouped selection ability of the penalties.
 They can be divided into two types. The first type contains light tailed errors models (1) whereas the second type is composed of heavy tailed errors models (2 and 3). Models 1 allows to quantify the deterioration of the performances of the robust algorithms  in the absence of outliers. Thinking about the maximum likelihood approach, the least squares loss (resp. Huber's loss) is well designed for Models 1 (resp. 2,3). 

\subsubsection{Assessing prediction methods}
\label{pred}

\begin{figure}[t!]
\vspace*{-6pt}
\centering
\hspace*{-1.8cm}\includegraphics{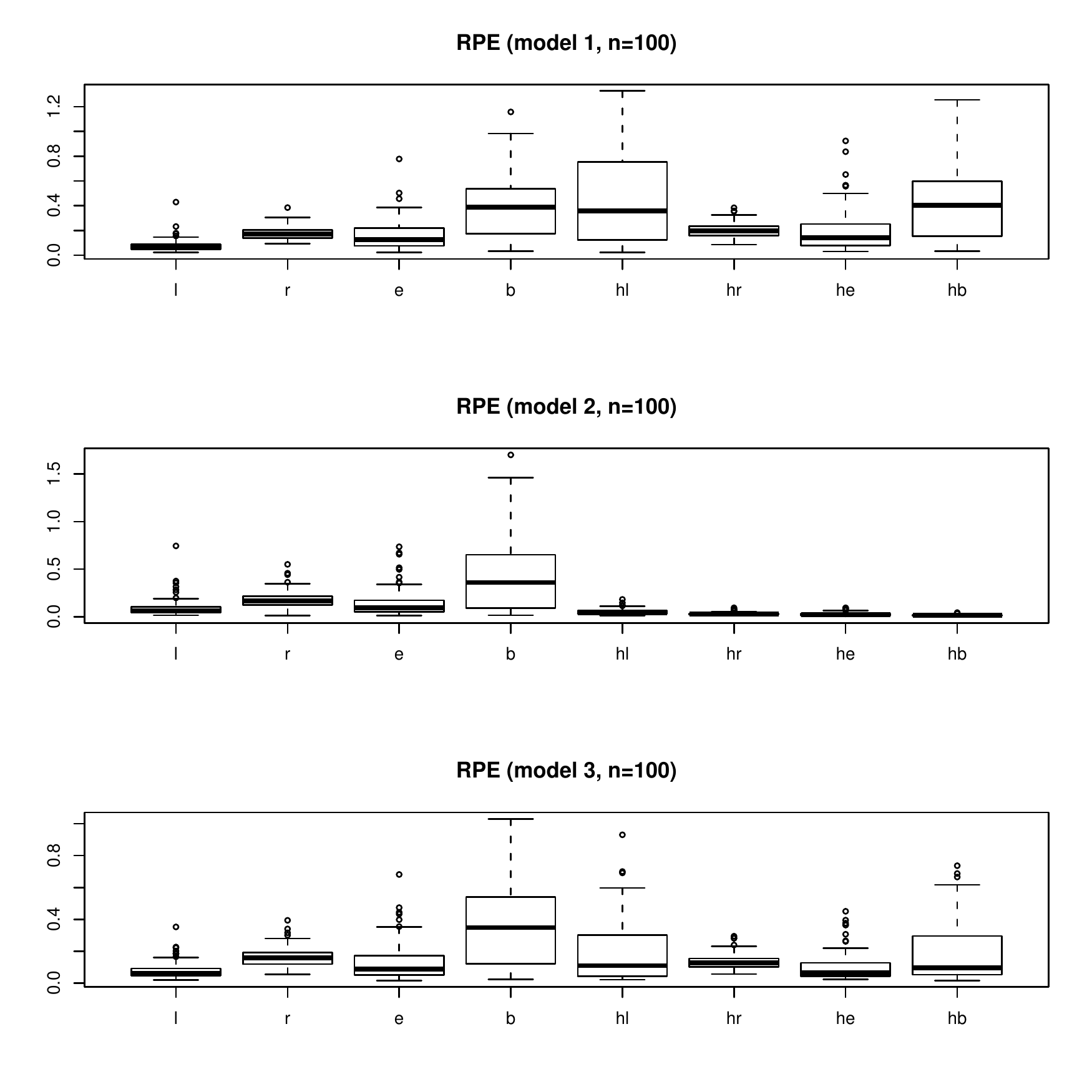}
\vspace*{-26pt}
\caption{For $n=100$, RPE for {\tt ad-lasso} (l), {\tt ridge} (r), {\tt ad-en} (e), {\tt ad-Berhu} (b), {\tt Huber-ad-lasso} (hl), {\tt Huber-ridge} (hr), {\tt Huber-ad-en} (he), and {\tt Huber-ad-Berhu} (hb). The boxplots are obtained without extreme values given by, for model 1 hl: 2.87; model 2 b: 2.95, hl: 2.94, he: 794.15; model 3 hl: 2.58.}
\label{RPE100}
\vspace*{-6pt}
\end{figure}

\begin{figure}[t!]
\vspace*{-6pt}
\centering
\hspace*{-1.8cm}\includegraphics{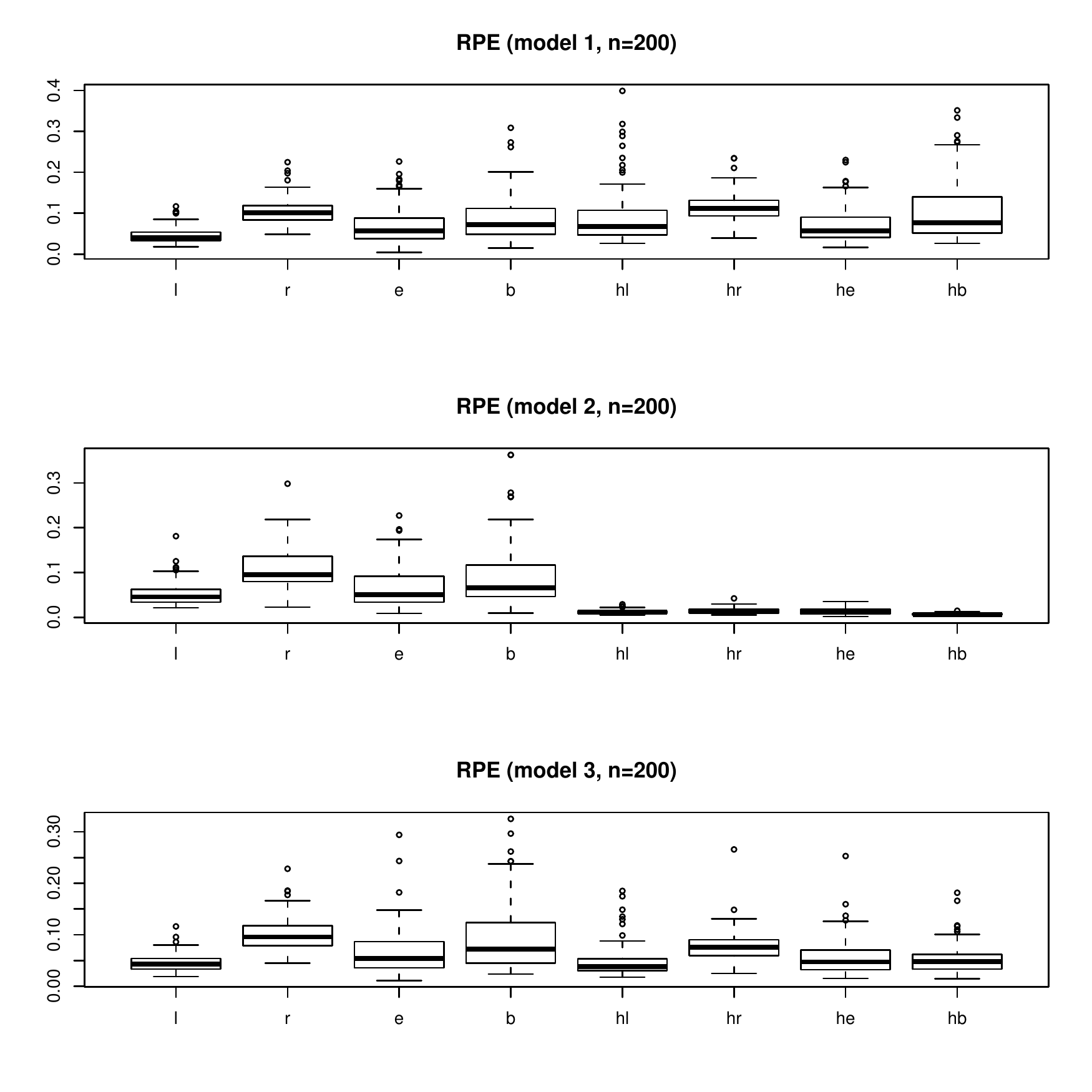}
\vspace*{-26pt}
\caption{For $n=200$, RPE for {\tt ad-lasso} (l), {\tt ridge} (r), {\tt ad-en} (e), {\tt ad-Berhu} (b), {\tt Huber-ad-lasso} (hl), {\tt Huber-ridge} (hr), {\tt Huber-ad-en} (he), and {\tt Huber-ad-Berhu} (hb). The boxplots are obtained without extreme values given by, for model 1 b: 2.95, hl: 2.48, 2.95,12.79,  2.86,  2.54,  2.96, 2.95; model 3 b: 2.95, 2.95, 2.95, hl: 2.95, 2.94, 2.95, 2.51, 49.03.}
\label{RPE200}
\vspace*{-6pt}
\end{figure}

\begin{figure}[t!]
\vspace*{-6pt}
\centering
\hspace*{-1.8cm}\includegraphics{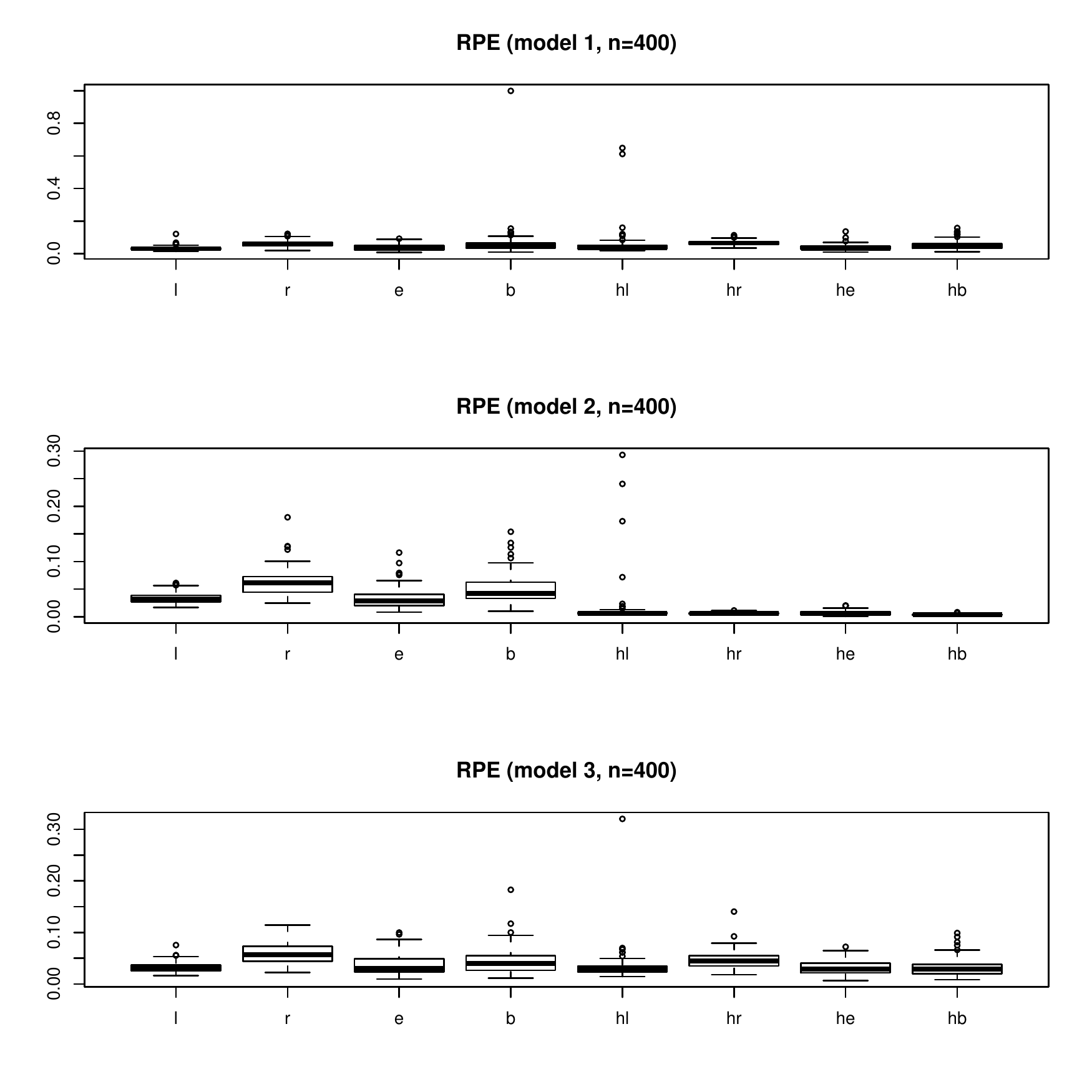}
\vspace*{-26pt}
\caption{For $n=400$, RPE for {\tt ad-lasso} (l), {\tt ridge} (r), {\tt ad-en} (e), {\tt ad-Berhu} (b), {\tt Huber-ad-lasso} (hl), {\tt Huber-ridge} (hr), {\tt Huber-ad-en} (he), and {\tt Huber-ad-Berhu} (hb). The boxplots are obtained without extreme values given by, for model 1 b: 2.95, 2.95, hl: 2.95, 2.49, 2.90, 2.95, 2.94, 2.95, 2.95, 2.93; model 2 b: 2.95, 2.95, 0.99; model 3 b: 2.95, 2.95, hl: 8.97, 2.95.}
\label{RPE400}
\vspace*{-6pt}
\end{figure}

To compare the performances of the various algorithms in the fixed design setting, the performances are measured both by the prediction errors and the model selection ability. For any considered underlying models,  we generate a first set of $n$ training designs ($\mathbf{x}_1,\cdots,\mathbf{x}_n$) and a second set of $m=$10 000 test designs ($\mathbf{x}_{n+1},\cdots,\mathbf{x}_{n+m}$).  These two sets are centered in mean to stick on the theoretical definition \eqref{eq:model} of the model (i.e. ensures that $\sum_{i=1}^n \mathbf{x}_i=0$). Since the theoretical results are established in fix design framework, the training and test design are fixed once and for all:  they  will be used for {\it all} the data generations. $100$ training sets of size $n$ are generated according to definition \eqref{eq:model} of the model. All the algorithms have been run on the $100$ training sets of size $n=$100, 200, 400 and their prediction capacity have been evaluated on the test design set of size $m=$10 000. To compare the prediction accuracy, the Relative Prediction Errors (RPEs) already considered in \cite{Zou06} are computed (see also \cite{Lam10} for explicit definition). Figures \ref{RPE100}, \ref{RPE200} and \ref{RPE400},  provide the boxplots associated with the 100 obtained RPE. 

\begin{figure}[t!]
\vspace*{-6pt}
\centering
\hspace*{-1.8cm}\includegraphics{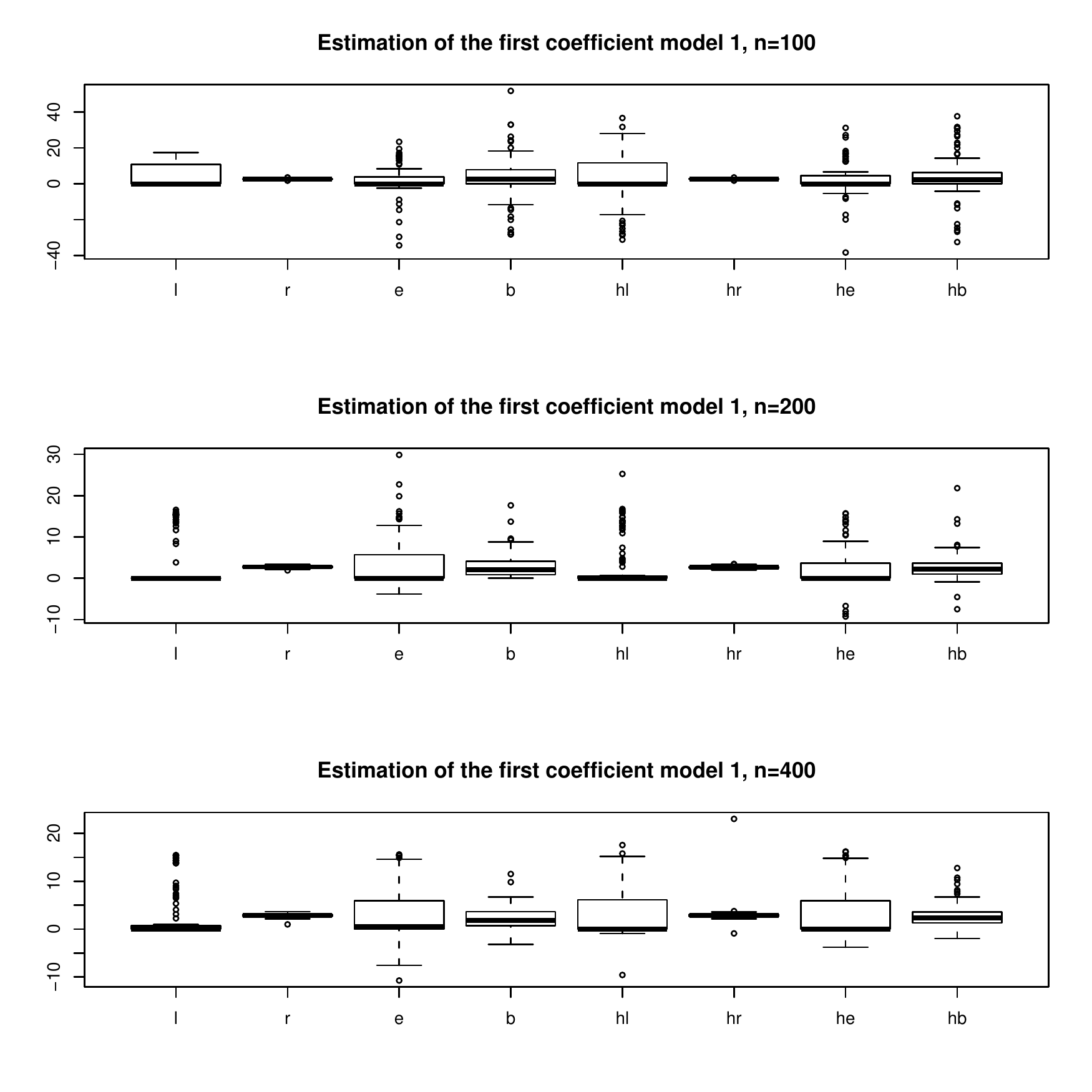}
\vspace*{-26pt}
\caption{Model 1: 
estimations of first influencing coefficient (true value is equal to 3) by {\tt ad-lasso} (l), {\tt ridge} (r), {\tt ad-en} (e), {\tt ad-Berhu} (b), {\tt Huber-ad-lasso} (hl), {\tt Huber-ridge} (hr), {\tt Huber-ad-en} (he), and {\tt Huber-ad-Berhu} (hb).}
\label{model1}
\vspace*{-6pt}
\end{figure}

\begin{figure}[t!]
\vspace*{-6pt}
\centering
\hspace*{-1.8cm}\includegraphics{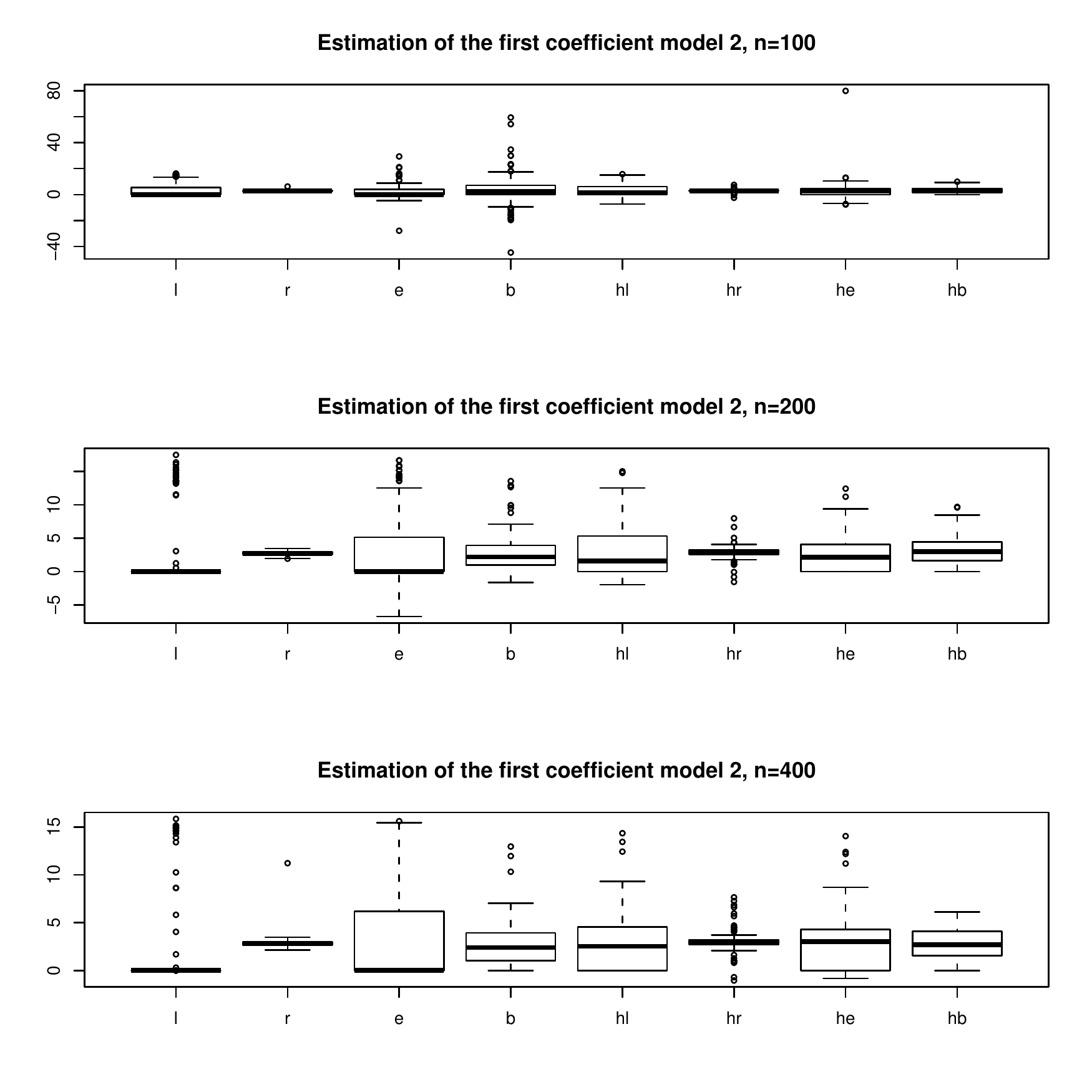}
\vspace*{-26pt}
\caption{Model 2: 
estimations of first influencing coefficient (true value is equal to 3) by {\tt ad-lasso} (l), {\tt ridge} (r), {\tt ad-en} (e), {\tt ad-Berhu} (b), {\tt Huber-ad-lasso} (hl), {\tt Huber-ridge} (hr), {\tt Huber-ad-en} (he), and {\tt Huber-ad-Berhu} (hb). }
\label{model2}
\vspace*{-6pt}
\end{figure}

\begin{figure}[t!]
\vspace*{-6pt}
\centering
\hspace*{-1.8cm}\includegraphics{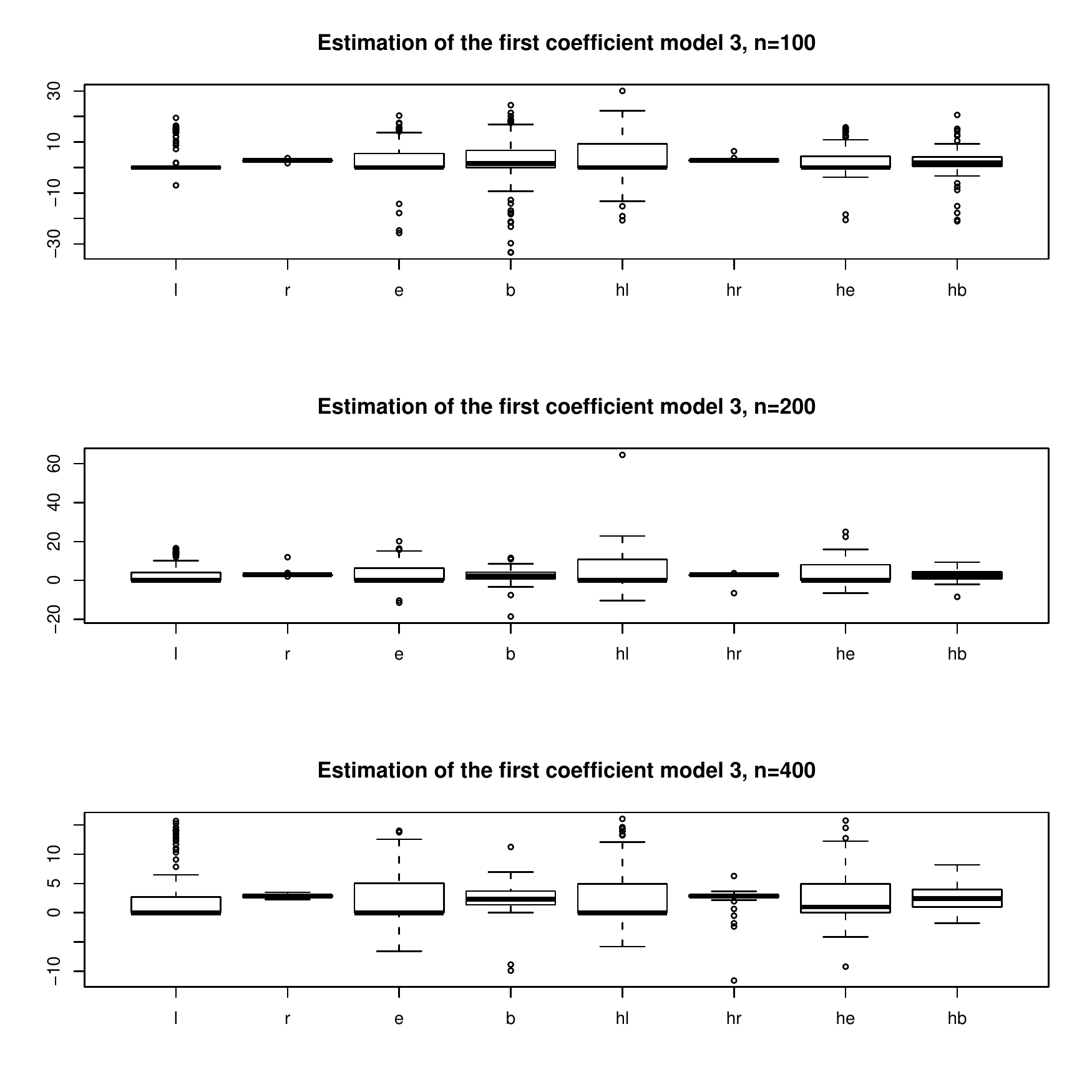}
\vspace*{-26pt}
\caption{Model 3: 
estimations of first influencing coefficient (true value is equal to 3) by {\tt ad-lasso} (l), {\tt ridge} (r), {\tt ad-en} (e), {\tt ad-Berhu} (b), {\tt Huber-ad-lasso} (hl), {\tt Huber-ridge} (hr), {\tt Huber-ad-en} (he), and {\tt Huber-ad-Berhu} (hb).}
\label{model3}
\vspace*{-6pt}
\end{figure}

The model selection ability of the algorithms are reported  in the same manner as done by  \cite{Wan07},  \cite{Tib96} and \cite{Fan01} in Tables \ref{modselecMod1}, \ref{modselecMod2} and \ref{modselecMod3}. Ridge penalty procedures are not reported since they do not constitute variables selection procedures.
To provide the indicators defined below, a coefficient is considered to be zero if it absolute value is strictly less than $10^{-5}$ (i.e. its five first decimals vanish).
In all cases, amongst the $100$ obtained estimators, the first column (C) counts the number of well chosen models i.e. the cases where  $15$ first coordinates of $\hat{\beta}$ are non-zeros {and} the $25$ last coefficients are zeros. To go further in the model selection ability analysis, we consider other measurements. The first  (in the second column (O)) represents the number of overfitting models (i.e. those selecting all the non-zeros coefficients {and} at least one zero coefficient). The second (in the third column (U)) reports the number of chosen underfitting models (i.e. those not selecting at least one non-zero coefficient). In this way, all the $100$ models are counted one time. Columns (0) and (U) aim to explain the results obtained in (C).  The column (Z) is the average number of estimated zeros, the column (CZ) provides the average number of correctly estimated zeros and (TZ) recall the theoretical zeros number. The column (CNZ) is  the average number of correctly estimated non zeros and (TNZ) recall the theoretical non zeros number. 
Models selection abilities are closely related to the accuracy of estimations of the coefficients. This fact is illustrated by boxplots  of the coefficients estimations (see Figures \ref{model1}, \ref{model2} and \ref{model3}).

Concerning the hyperparameter choices, the regularization parameters associated with adaptive lasso or Berhu penalties are chosen by BIC criterion on {each} of the $100$ training sets as described at Section~\ref{ss:hyperparameter}. The same grid has always been used for each method. It is composed of 100 points log-linearly spaced between 0 and 1400 for Berhu and 200  points log-linearly spaced between 0 and 10 000 for lasso. For Huber's loss, the simulation studies report the performances obtained with $M=1.345$. This value has been recommended by Huber in \cite{Hub81}. For adaptive Berhu penalty, we report the performances obtained with $L=M=1.345$
Let us remark that it is possible to chose the $M$ and $L$ parameters from the data (for example by cross-validation simultaneous with the tuning parameter). But in practice we do not observe some improvement to make it data adaptive. For ridge-type procedures, the hyperparameter is chosen as usually by 5-fold cross-validation on each of the 100 training sets. The grid is composed of 100 points log-linearly spaced between 0 and 1400. For en-type procedure, we use the similar protocol as in \cite{Zou05}: we first pick a relatively small grid of values for $\lambda_{2,n}$ over $\{0$, 0.01, 0.1, 1, 10, 100$\}$ and 25 points log-linearly spaced between 0 and 5000 for $\lambda_{1,n}$. Then the both parameters are 
chosen simultaneously by 5-fold cross-validation.

\subsubsection{Comparison results}
\label{comp}

Tables \ref{modselecMod1}, \ref{modselecMod2} and \ref{modselecMod3} present the performances in terms of selection model ability. First we see that whatever the model the behavior of the methods are the same.
The lasso and en penalties methods lead in general to underfitting models (columns U). It is surprising for the en penalty. Indeed the penalty imposed on the small coefficients is the sum of an $\ell_1$-norm and a squared $\ell_2$-norm. 
This implies that the obtained penalty is closer to differentiability than the $\ell_1$-penalty. As shown in \cite{AntoFan2001}, if the penalty is far from differentiability, more small  coefficients  are deleted. For these examples, the en penalty as the same behavior as the lasso one. As a consequence, these penalties have a relatively high number of zeros with correct zeros number (columns Z) very close to the true one (columns TZ). But the correct non zeros number (columns CNZ) is very low in comparison with the true one (columns TNZ). The fact that en and  lasso type methods underfit is reduced for Model 2 and for Huber loss. In all cases, these methods almost never identify the right model.
The Behru penalty leads to some compromise between over and under fitting. We point out that contrary to en and lasso type methods, there is a case where Berhu type method identifies the right model a reasonable number of times: it is  Model 2 with Huber loss. It is a little less good in terms of correct zeros but much better in terms of non zeros number. 

This behavior occurs on the quality of estimation of the non zero coefficients (see Figures \ref{model1}, \ref{model2} and \ref{model3})). Let us note that we have only considered the first coefficient $\beta_1$ and that the conclusions for the other non zero coefficient are the same. The ridge method is given here as a reference since it is known to lead good performances in presence of high correlation between the covariables. We observe that the Berhu penalty lead to good performance in terms of bias as ridge with higher variability than the ridge one. The bias and sometimes the variability are very high for the other methods du to their tendency to underfitting.

Figures \ref{RPE100}, \ref{RPE200} and \ref{RPE400},  provide the boxplots associated with RPE. 
As expected, the ordinary least squares loss leads to better performance for Model 1 (excepted for $n=400$) and leads to less good performance for the Model 3 and especially for the Model 2. We observe that 
{\tt ad-Berhu} and {\tt Huber-ad-lasso} provide several extreme values du to numerical instabilities and are often more variable. 

\subsection{Prostate cancer data example}

This data set comes from a prostate cancer study (see \cite{Sta89}) and analyzed earlier in the elastic net paper by \cite{Zou05,Gosh2007}. There are eight clinical covariates namely: logarithm of the cancer volume (lcavol), logarithm of the prostate weight (lweight), age, the logarithm of the amount of benign prostatic hyperplasia (lbph), seminal vesicle invasion (svi), logarithm of the capsular penetration (lcp), Gleason score (gleason) and percentage Gleason score 4 or 5 (pgg45). The response is the logarithm of prostate-specific antigen (lpsa). The predictors are are named as $1,\ldots ,8$ in results. OLS and the previous methods were applied to these data. 

\begin{figure}[t!]
\vspace*{-6pt}
\centering
\hspace*{-1.8cm}\includegraphics{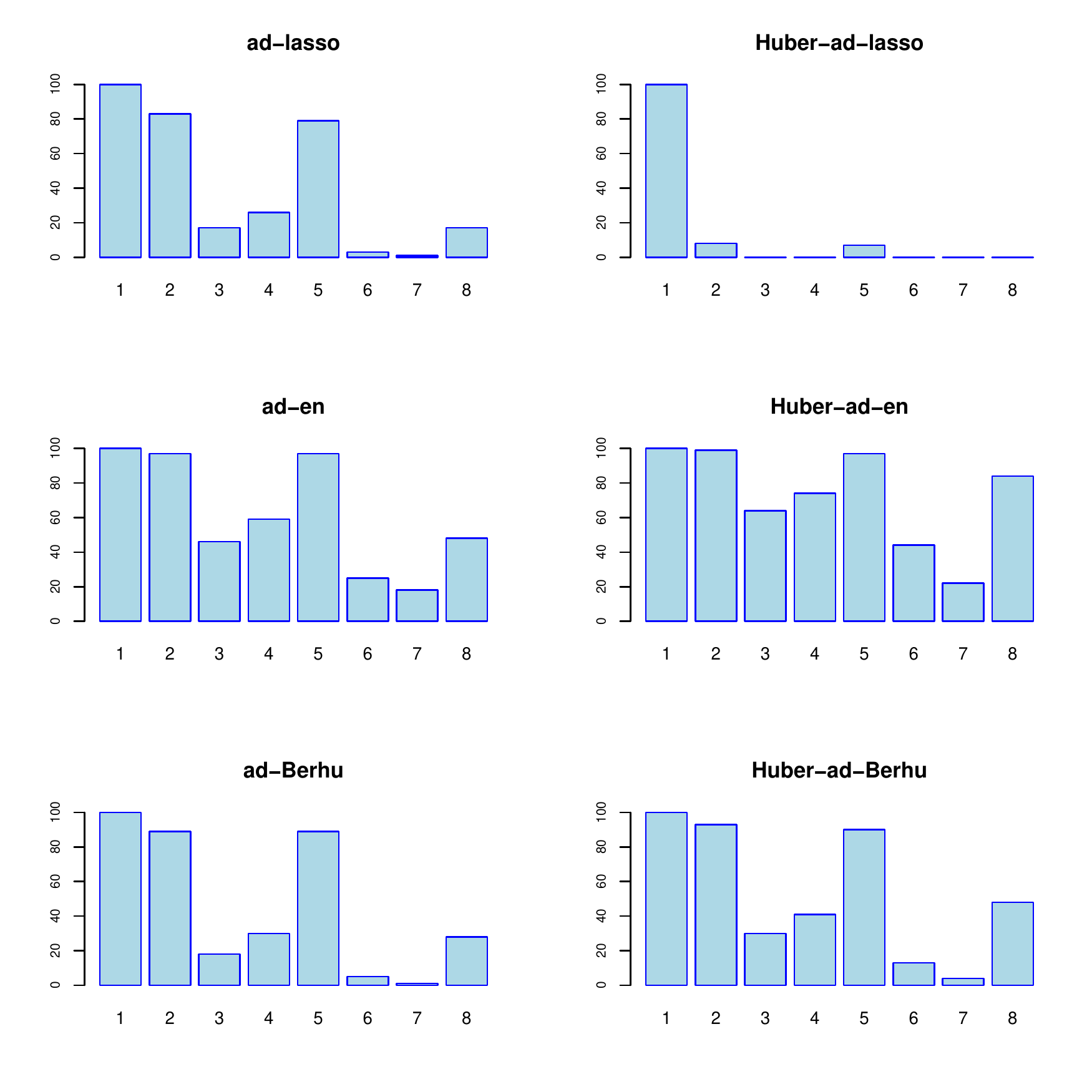}
\vspace*{-26pt}
\caption{Prostate cancer data: histogram associated with number of selection of each variables in the re-sampling study.}
\label{histo}
\vspace*{-6pt}
\end{figure}

In \cite{Zou05}, the data were divided into to parts: a training set with 67 observations and a test set with 30 observations while in \cite{Gosh2007}, they have divided (randomly) the original data set in to training and testing set containing 60 and 37 observations respectively. To fairly compare the methods we propose to perform a resampling study: we have divided $100$ times (randomly) the original data set into training and testing set containing 67 and 30 observations respectively. The hyperparameters are chosen as in the simulation study. We then compared the performances of the methods by computing their RPE on the 100 resampling testing sets (see Table~\ref{prederr}).  	
Contrary to what had been observed in \cite{Zou05,Gosh2007}, our resampling study does not allow us to claim that one method emerges in terms of RPE: almost all these methods have similar RPE. 
We can only say perhaps  {\tt Huber-ad-lasso}  is slightly less good. Let us notice that we observe a great variability in the choice of $\lambda_{2,n}$  for the adaptive en-type procedures (see first column of Table~\ref{prederr}). This is  also the case for {\tt Huber-ad-lasso}. As a contrary, the choice of $\lambda_{n}$ for Behru type procedures is more stable (it is comparable to the stability of ridge).
Figure~\ref{histo} show (except for OLS and ridge procedures) the histogram associated with the selected variables. We see that Berhu penalties leads to good models in terms of sparsity in comparison with en penalties.  We observe that Behru type procedures are compromise between lasso type methods which select too few variables and en type methods which select too many variables.

\section{Appendix}

\subsection{Computations: software used for numerical optimization}
\label{ap:algo}

When the regularization parameter is fixed, to solve all the involved optimization problems we used \texttt{CVX}, a package for specifying and solving convex programs \cite{CVX,DCP}. \texttt{CVX} is a set of Matlab functions using the methodology of disciplined convex programming. Disciplined convex programming imposes a limited set of conventions or rules,
which are called the DCP ruleset. Problems which adhere to the ruleset can be rapidly
and automatically verified as convex and converted to solvable form. Problems that
violate the ruleset are rejected, even when convexity of the problem is obvious to the
user. The version of \texttt{CVX} we use, is a preprocessor for  the convex optimization solver SeDuMi (Self-Dual-Minimization \cite{SeDuMiSturm1999}).

Let us now recall a well-known fact of convex analysis: the Huber function is the Moreau-Yosida regularization of the absolute value function (\cite{HiLe91,Rock70,Sard01}). Precisely, it can be easily shown that the Huber function satisfies
$$
{\cal H}_M(z)=\min_{v\in\rset}\left((z-v)^2+2M\vert v\vert\right)\,.
$$
We can derive the same kind of formulation for the BerHu function leading to a characterization of the BerHu function as quadratic optimization problem. Indeed, the function \eqref{defBerHu} satisfies 
$$
{\cal B}_L(z)=\min_{w\geq L\vee\vert z\vert}\left(\frac{w^2}{2L}-w+\vert z\vert+\frac{L}{2}\right)\,,
$$
where $a\vee b$ denotes the maximum of the two real numbers $a$ and $b$. The proof of this equality is trivial since it amounts to minimize a quadratic function on an interval.

This allows to write our optimization problem in a conforming manner to use \texttt{CVX}. Note that \cite{Owe06} uses an expression of ${\cal H}_M(z)$ as the solution of a quadratic optimization problem  (borrowed from the user guide of \texttt{CVX}) to  write his problem in a conforming manner to use \texttt{CVX}. However, the expression of \cite{Owe06} involves more constraints and more variables than the previous formulation.
We give here the way to use \texttt{CVX} in order to compute the estimators {\tt{alpha}}=$\hat{\alpha}^{{\cal H}adl}$, 
{\tt{beta}}=$\hat{\beta}^{{\cal H}adl}$ and {\tt{s}}=$\hat{s}^{{\cal H}adl}$. The variable {\tt{X}} represents the design matrix $\X$.
The unpenalized estimator {\tt{betaUNP}}$=\hat{\beta}_{{\cal H}}$ is calculated beforehand (using also \texttt{CVX}) and the regularisation parameter $\lambda_n$ is fixed and denoted by {\tt{lambda}}.
\begin{verbatim}
cvx_begin
variables alpha beta(p) s v(n) tau w(p);
minimize (n*s+quad_over_lin(y-alpha-X*beta-v,s)+2*M*norm(v,1)
+ mu*(tau*norm(betaUNP,1)+quad_over_lin(w./(sqrt(abs(betaUNP))),2*L*tau)
+norm(beta./betaUNP,1)-sum(w./abs(betaUNP))+0.5*L*tau*norm(1./betaUNP,1)))
subject to
s > 0;
tau > 0;
w >= L*tau;
w >= abs(beta);
cvx_end
\end{verbatim}
Let us remark that {\tt{betaUNP}} is computed in the same way but deleting the term multiplied by {\tt{lambda}}.

\subsection{Proof of Theorem~\ref{theo:groupingeffect}}
Since $\hat{\beta}^{adb}_i\neq0$, we have $\hat{\beta}^{adb}\neq0$ and $\hat{\tau}^{adb}>0$.
Consequently, the definition of partial derivatives involving Newton's quotient leads to the following KKT conditions by differentiating with respect to  $\beta_i$, $\beta_j$ and $\tau$ :
\begin{equation}
 \label{i}
-2{x_{i}}^T\left(\underline{y}-\hat{\alpha}^{adb}\indic{n}-\X\hat{\beta}^{adb}\right)+{\lambda_n}\hat{w}^{adb}_i{\cal B}_L'\left(\frac{\hat{\beta}^{adb}_i}{\hat{\tau}^{adb}}\right)=0,
\end{equation}
\begin{equation}
 \label{j}
-2{x_{j}}^T\left(\underline{y}-\hat{\alpha}^{adb}\indic{n}-\X\hat{\beta}^{adb}\right)+{\lambda_n}\hat{w}^{adb}_j{\cal B}_L'\left(\frac{\hat{\beta}^{adb}_j}{\hat{\tau}^{adb}}\right)=0,
\end{equation}
$$
\sum_{j:\vert\hat{\beta}^{adb}_j\vert>L\hat{\tau}^{adb}}\hat{w}^{adb}_j\left(
\frac 1 {2L} \left(\frac{\hat{\beta}^{adb}_j}{\hat{\tau}^{adb}}\right)^2-\frac{L}{2}\right)=\sum_{j=1}^p \frac 1{\hat{w}^{adb}_j}.
$$
The last score equation implies that the set $G=\{j\in [1;p], \vert\hat{\beta}^{adb}_j\vert>L\hat{\tau}^{adb}\}$ is non-empty. Let us now distinguish some cases involving this set on indices. To begin with, if both the indexes $i$ and $j$ belong to $G$, equations \eqref{i} and  \eqref{j}  become 
\begin{equation}
\label{zz}
-2{x_{i}}^T\left(\underline{y}-\hat{\alpha}^{adb}\indic{n}-\X\hat{\beta}^{adb}\right)+
{\lambda_n}\hat{w}^{adb}_j\frac{\hat{\beta}^{adb}_i}{L\hat{\tau}^{adb}}=0,
\end{equation}
and
$$
-2{x_{j}}^T\left(\underline{y}-\hat{\alpha}^{adb}\indic{n}-\X\hat{\beta}^{adb}\right)+{\lambda_n}\hat{w}^{adb}_j\frac{\hat{\beta}^{adb}_i}{L\hat{\tau}^{adb}}=0.
$$
Substracting the second one to the first one and using Cauchy-Schwarz inequality, we get :
$$
\vert\hat{w}^{adb}_i\hat{\beta}^{adb}_i-\hat{w}^{adb}_j\hat{\beta}^{adb}_j\vert\leq 
\frac{2L\hat{\tau}^{adb}}{\lambda_n}\left\|x_i-x_j\right\|_2\left\|\underline{y}-\hat{\alpha}^{adb}\indic{n}-\X\hat{\beta}^{adb}\right\|_2\,.
$$
The definition of  ($\hat{\alpha}^{adb},\hat{\beta}^{adb},\hat{\tau}^{adb}$)  as a minimizer implies that, for all $\tau>0$,
$$
\left\|\underline{y}-\hat{\alpha}^{adb}\indic{n}-\X\hat{\beta}^{adb}\right\|_2\leq \left\| \underline{y}\right\|_2+
\lambda_n \tau\sum_{j=1}^{p}\frac 1 {\hat{w}^{adb}_j}.
$$
Now, letting $\tau$ tends to $0$ in this inequality, we get:
\begin{equation}
\label{y}
\left\|\underline{y}-\hat{\alpha}^{adb}\indic{n}-\X\hat{\beta}^{adb}\right\|_2\leq\left\| \underline{y}\right\|_2.
\end{equation}
This leads to equation \eqref{to} of the Theorem~\ref{theo:groupingeffect} since $C_{i,j}=1$ in this case.

Next, let us consider the case where only one index among $\{i,j\}$ belongs to G. 
If $i$ and $j$ are switched (if necessary), we can suppose that $i\in G$ and $j\notin G$. In this case,  equations \eqref{i} and   \eqref{j}  become \eqref{zz} and 
$$
-2{x_{j}}^T\left(\underline{y}-\hat{\alpha}^{adb}\indic{n}-\X\hat{\beta}^{adb}\right)+
{\lambda_n}\hat{w}^{adb}_j {\rm sign}\left(\hat{\beta}^{adb}_i\right)=0.
$$
These two equalities lead to 
$$
\hat{w}^{adb}_i\hat{\beta}^{adb}_i-\hat{w}^{adb}_j\hat{\beta}^{adb}_j
=
\frac{2L\hat{\tau}^{adb}}{\lambda-n}\left(x_i-\frac{\vert\hat{\beta}^{adb}_j\vert}{L\hat{\tau}^{adb}}x_j\right)^T
\left(\underline{y}-\hat{\alpha}^{adb}\indic{n}-\X\hat{\beta}^{adb}\right)\,.
$$
Combining Cauchy-Schwarz inequality and inequality \eqref{y}, this leads to
$$
\vert\hat{w}^{adb}_i\hat{\beta}^{adb}_i-\hat{w}^{adb}_j\hat{\beta}^{adb}_j\vert\leq 
\frac{2L\hat{\tau}^{adb}}{\lambda_n}\vert\vert\underline{y}\vert\vert_2\sqrt{\vert\vert x_i\vert\vert_2^2+\vert\vert x_j\vert\vert_2^2-2\frac{\vert\hat{\beta}^{adb}_j\vert}{L\hat{\tau}^{adb}}{x_i}^Tx_j}\,,
$$
where we have used $j\notin G$. This implies equation \eqref{to} of the Theorem~\ref{theo:groupingeffect}  since $C_{i,j}={\vert\hat{\beta}^{adb}_j\vert}/{L\hat{\tau}^{adb}}$ in this case.

Finally, when $i$ and $j$ do not belong to $G$, using similar arguments we obtain
$$
\vert\hat{w}^{adb}_i\hat{\beta}^{adb}_i-\hat{w}^{adb}_j\hat{\beta}^{adb}_j\vert\leq 
\frac{2L\hat{\tau}^{adb}}{\lambda_n}\vert\vert\underline{y}\vert\vert_2\sqrt{\vert\vert x_i\vert\vert_2^2+\vert\vert x_j\vert\vert_2^2-2\frac{\vert\hat{\beta}^{adb}_j\hat{\beta}^{adb}_i\vert}{L^2\hat{\tau}^{adb2}}{x_i}^Tx_j}\,,
$$
that implies equation \eqref{to} of the Theorem \eqref{to}  since $C_{i,j}={\vert\hat{\beta}^{adb}_j\hat{\beta}^{adb}_i\vert}/({L^2\hat{\tau}^{adb2}})$ in this case. $\blacksquare$

\subsection{Proof of Theorem~\ref{8}}\label{approofthe2}

The asymptotic normality of this estimator is proved in Step 1 and
the consistency in variable selection in the Step 2. This proof
in an adaptation to our case of the proof given by  \cite{Zou06} or \cite{Lam10}. 
The difference with \cite{Lam10} concerns the treatment of the penalty term. So in the
following, we will use notations similar to the ones of \cite{Lam10}.
We will point out the difference between the both proofs.

{\bf Step 1}. Let us first prove the asymptotic normality. 
Let us define $U_n(u)= Q^{{\cal H}adb}\left((\alpha^*,\beta^*,s^*,\tau^*\right
)^T+u/\sqrt{n})-Q^{{\cal H}adb}(\alpha^*,\beta^*,s^*,\tau^*)$ with
$u=(u_0,\ldots,u_{p+2})^T\in\rset^{p+3}.$
Obviously, $U_n(u)$ is minimized at
\[
\hat{u}^{(n)}=\sqrt{n}\left(\hat{\alpha}^{{\cal H}adb}-\alpha^*,\hat
{\beta}^{{\cal H}adb}-\beta^*,\hat{s}^{{\cal H}adb}-s^*, \frac {\sqrt{\lambda_n}}{\sqrt{n}}(\hat{\tau}^{{\cal H}adb}-\tau^*)\right)^T.
\]
The principle of the proof of \cite{Zou06} or \cite{Lam10} is to study the epi-limit
of $U_n$. Using the proof of theorem 3.2 in \cite{Lam10},  we only need to study the epi-limit of the penalty 
term given by
$$P_n(u)=\lambda_n\left(\widetilde{P}^{adb}\left(\beta^*+\frac{u_{1:p}}{\sqrt{n}},\tau^*+\frac{u_{p+2}}{\sqrt{\lambda_n}}\right)-\widetilde{P}^{adb}(\beta^*,\tau^*)\right),$$
where $\widetilde{P}^{adb}(\beta,\tau)={P}^{adb}(\beta,\tau),$ if $\tau\geq 0,$ $\infty$ if $\tau<0.$ The epi-limit of this term is given in the Lemma~\ref{cvpen4} .
This lemma together with lemma~2 of \cite{Lam10} indicates that $U_n\rightarrow_{e-d}U,$ where
$U(u)=A_{s^*}\left(u_{1:p}^TVu_{1:p}+u_0^2\right
)+D_{s^*}{u^2_{p+1}}-W^Tu+{u_{p+2}}^2C(u_{p+2}),$ if $u_j=0,$ $\forall\,j\notin\CA$, $+\infty$ otherwise. Under condition $\beta^*\neq0,$ equation \eqref{proptaustar} in Lemma~\ref{existtaustar} implies that $\sum_{j=1}^{p_0}\vert\beta^*_j\vert^{2-\gamma}\indic{\vert\beta^*_j\vert>L\tau*}>0
$ thus the function $z\rightarrow {z}^2C(z)$ is strictly convex. Moreover, $V_{1,1}$ is supposed  positive definite in assumption {\bf (D2)} and we assume that the noise satisfies {\bf(N2)}. Consequently, $U$ get a unique argmin and  the asymptotic normality part is proved.

{\bf Step 2.} Let us now show the consistency in variable selection part.  It suffices to show that $\Pb{\CA\subset\CA_n}\rightarrow1$ as $n$  tends to infinity and $\Pb{\CA^c\subset{\CA_n}^c}\rightarrow1$ as $n$  tends to infinity. The first claim is an easy consequence of asymptotical normality obatined in Step~1. 

Let us now show the second claim.  Let $j$ such that $\beta^*_j=0$. We have to prove that $\Pb{\hat{\beta}^{{\cal H}adb}_{j}\neq0}\rightarrow0$ as $n$ tends to infinity. As in \cite{Lam10}, we have for a such $j,$ 
\begin{multline*}
\Pb{\hat{\beta}^{{\cal H}adb}_{j}\neq0}\leq \Pb{\left(\hat{s}^{{\cal H}adb},\hat{\tau}^{{\cal H}adb}\right)=(0,0)}
+\\\Pb{\hat{\tau}^{{\cal H}adb}>0\,\text{and}\,\hat{s}^{{\cal H}adb}>0\,
\text{and}\,\sum_{i=1}^n{x}_{i,j}\dHu{\frac{y_i-\hat{\alpha}^{{\cal H}adb}-\mathbf{x}_i^T\hat{\beta}^{{\cal H}adb}}{\hat{s}^{{\cal H}adb}}}{M}=-\lambda_n\hat{{w}}^{{\cal H}adb}_j\dIH{\frac{\hat{\beta}_{j}^{{\cal H}adb}}{\hat{\tau}^{{\cal H}adb}}}{L}}.
\end{multline*}
Using similar arguments as in \cite{Lam10}, we have, as $n$ tends to infinity, 
$$\Pb{(\hat{s}^{{\cal H}adb},\hat{\tau}^{{\cal H}adb})=(0,0)}\rightarrow0.$$
Since $\forall x\in\rset^*$,\, $\vert\dIH{x}{L}\vert\geq1$, we have
\begin{eqnarray}
\nonumber \Pb{\hat{\tau}^{{\cal H}adb}>0\,\text{and}\,\hat{s}^{{\cal H}adb}>0\,\text{and}\,\sum_{i=1}^n\mathbf{x}_{i,j}\dHu{\frac{y_i-\hat{\alpha}^{{\cal H}adb}-\mathbf{x}_i^T\hat{\beta}^{{\cal H}adb}}{\hat{s}^{{\cal H}adb}}}{M}=-\lambda_n\hat{{w}}^{{\cal H}adb}_j\dIH{\frac{\hat{\beta}_{j}^{{\cal H}adb}}{\hat{\tau}^{{\cal H}adb}}}{L}}\\\nonumber\leq \Pb{\hat{s}^{{\cal H}adb}>0\,\text{and}\,\frac{1}{\sqrt{n}}\left |\sum_{i=1}^n {x}_{i,j}\dHu{\frac{y_i-\hat{\alpha}^{{\cal H}adb}-\mathbf{x}_i^T\hat{\beta}^{{\cal H}adb}}{\hat{s}^{{\cal H}adb}}}{M}\right |\geq\frac{\lambda_n}{\sqrt{n}}\hat{{w}}^{{\cal H}adb}_j}
\end{eqnarray}
As in \cite{Lam10}, we have 
$$\frac{1}{\sqrt{n}}\sum_{i=1}^n{x}_{i,j}\dHu{\frac{y_i-\hat{\alpha}^{{\cal H}adb}-\mathbf{x}_i^T\hat{\beta}^{{\cal H}adb}}{\hat{s}^{{\cal H}adb}}}{M}=O_P(1),$$ 
and ${\sqrt{n}}/({\lambda_n\hat{{w}}^{{\cal H}adb}_j})\stackrel{\mathbb{P}}{\rightarrow}0,$ that 
implies that  $\Pb{\hat{\beta}^{{\cal H}adb}_{j}\neq0}\rightarrow0$ as $n$ tends to infinity.
$
\blacksquare
$

\subsection{Technical lemma}
\subsubsection{Proof of lemma~\ref{cvpen4}}\label{ss:lemma2}
\begin{lemme}
\label{cvpen4} Suppose that $\lambda_n/n^{\gamma\wedge 1/2}\rightarrow0,$  $\lambda
_n n^{(\gamma-1)/2}\rightarrow\infty,$ $\lambda
_n \rightarrow\infty,$ $\lambda_n>1/3$ and $\beta^*\neq 0.$ Then we have 
$$
P_n(u)\rightarrow_{e-d}
\left\{
\begin{array}{ll}
{u_{p+2}}^2C(u_{p+2})& \text{if }u_j=0,\,\forall\,j\notin\CA,\\
+\infty&\text{otherwise}\,,
\end{array}
\right.
$$
where 
$$C(u_{p+2})=\frac{1}{2L{\tau^*}^3}\sum_{j=1}^{p_0}\vert\beta^*_j\vert^{2-\gamma}\indic{\vert\beta^*_j\vert>L\tau*}+\frac{L^{1-\gamma}}{2{\tau^*}^{(\gamma+1)}}\#\{1\leq j\leq p,\,\vert\beta^*_j\vert=L\tau^*\}\indic{u_{p+2}<0}.$$
\end{lemme}

Since $\beta^*\neq0$, Lemma~\ref{existtaustar} ensures that $\tau^*>0$. Consequently, we have  
$P_n(u)=\sum_{j=1}^{p}P_{n,j}(u)$, where
$$
P_{n,j}(u)=\left\{
\begin{array}{lll}
\lambda_n\left(\frac{u_{p+2}}{\sqrt{\lambda_n}\hat{w}^{adb}_j}+ \hat{w}^{adb}_j\left(\tau^*+\frac{u_{p+2}}{\sqrt{\lambda_n}}\right){\cal B}_L\left(\frac {\beta^*_j+\frac{u_j}{\sqrt{n}}}{\tau^*+\frac{u_{p+2}}{\sqrt{\lambda_n}}}\right)-\tau^*\hat{w}^{adb}_j{\cal B}_L\left(\frac{\beta^*_j}{\tau^*}\right)\right)&\text{if}& u_{p+2}>-\sqrt{\lambda_n}\tau^*,\\
-\lambda_n\tau^*\left(\frac 1 {\hat{w}^{adb}_j}+\hat{w}^{adb}_j
{\cal B}_L\left(\frac{\beta^*_j}{\tau^*}\right)\right)&\text{if}&  u_{p+2}=-\sqrt{\lambda_n}\tau^*,\\
&\text{and}&u_j=-\sqrt{n}\beta^*_j,\\
+\infty&&\text{otherwise}.\\
\end{array}
\right.
$$
{\bf Step 1}. First let us prove that
\begin{equation}
\label{non08}
\sum_{j=1}^{p_0}P_{n,j}(u)\rightarrow_{e-d}u_{p+2}^2C(u_{p+2})\,.
\end{equation}
We show that, for every $u$ fixed in $\rset^{p+2}$, we have this convergence in probability.
Since $\tau^*>0$ and $\lambda_n\rightarrow+\infty$ as $n$ tends to infinity, for $n$ sufficiently large (with respect to a bound depending on  $u_{p+2}$), ${u_{p+2}}/{\sqrt{\lambda_n}}+\tau^*>0$ and 
$$
P_{n,j}(u)=
\frac {u_{p+2}\sqrt{\lambda_n}} {\hat{w}^{adb}_j} + {\lambda_n}{\hat{w}^{adb}_j}\left(G_j\left(\frac{u_j}{\sqrt{n}},\frac{u_{p+2}}{\sqrt{\lambda_n}}\right)-G_j(0)\right)\,,
$$
where 
$$
\forall\,j\in\,[1,p_0],\,G_j:(z_1,z_2)\rightarrow \left(z_2+\tau^*\right){\cal B}_L\left(\frac {z_1+\beta^*_j}{z_2+\tau^*}\right).
$$
For $1\leq j\leq p_0$ such that $\vert\beta^*_j\vert\neq L\tau^*$,  $G_j$ is two times differentiable at $0$ and the Taylor-Young  theorem entails that, $\forall\, (z_1,z_2)\,\in\,\rset^{2}$,
\begin{multline*}
G_j(z_1,z_2)=G_j(0)+z_1\dIH{\frac{\beta^*_j}{\tau^*}}{L}+z_2B\left(\frac{\beta^*_j}{\tau^*}\right)+\frac{z_1^2}{2L\tau^*}\indic{\vert\beta^*_j\vert>L\tau*}\\+\frac{z_2^2{\beta^*_j}^2}{2L{\tau^*}^3}\indic{\vert\beta^*_j\vert>L\tau*}-\frac{z_1z_2\beta^*_j}{L{\tau^*}^2}\indic{\vert\beta^*_j\vert>L\tau*}+\xi(z_1,z_2)\,,
\end{multline*}
where ${\xi(z_1,z_2)}/{\|(z_1,z_2\|^2}\rightarrow 0$ as $(z_1,z_2)\rightarrow 0$\,,  $B:z\in \rset \rightarrow \IH{z}{L}-z\dIH{z}{L}$ and we have used that $\ddIH{{\beta^*_j}/{\tau^*}}{L}=\indic{\vert\beta^*_j\vert>L\tau*}/L$. Consequently, for $1\leq j\leq p_0$ such that $\vert\beta^*_j\vert\neq L\tau^*$,
\begin{equation}
\label{nonL8}
P_{n,j}(u)=\frac{u_{p+2}\sqrt{\lambda_n}}{\hat{w}^{adb}_j}+\sqrt{\lambda_n}u_{p+2}{\hat{w}^{adb}_j}B\left(\frac{\beta^*_j}{\tau^*}\right)+\frac{u_{p+2}^2\vert\beta^*_j\vert^{2-\gamma}}{2L{\tau^*}^3}\indic{\vert\beta^*_j\vert>L\tau*}+a_{n,j}(u)\,,
\end{equation}
where 
$$a_{n,j}(u)=\frac{\lambda_n u_j \hat{w}^{adb}_j }{\sqrt{n}}\dIH{\frac{\beta^*_j}{\tau^*}}{L}+
\frac{\lambda_nu_j^2 \hat{w}^{adb}_j }{2nL\tau^*}\indic{\vert\beta^*_j\vert>L\tau*}-
\frac{\sqrt{\lambda_n}u_{p+2}\beta^*_ju_j \hat{w}^{adb}_j }{\sqrt{n}L{\tau^*}^2}\indic{\vert\beta^*_j\vert>L\tau*}
+\lambda_n\hat{w}^{adb}_j \xi\left(\frac{u_j}{\sqrt{n}},\frac{u_{p+2}}{\sqrt{\lambda_n}}\right).$$

Let us now consider $1\leq j\leq p_0$ such that $\vert\beta^*_j\vert= L\tau^*$. When $\beta^*_j= L\tau^*$, for $n$ sufficiently large (with respect to a bound depending on  $u$),
$$
P_{n,j}(u)=
\frac{\sqrt{\lambda_n}u_{p+2}} {\hat{w}^{adb}_j}
+{\lambda_n\hat{w}^{adb}_j}
\left(\left(\tau^*+\frac{u_{p+2}}{\sqrt{\lambda_n}}\right){\cal B}_L\left(L\tau^*+\frac {\frac{u_j}{\sqrt{n}}}{\tau^*+\frac{u_{p+2}}{\sqrt{\lambda_n}}}\right)-L\tau^*\right)
$$
Let us consider $n$ sufficiently large (with respect to a bound depending on $u$) such that $L\tau^*+{u_j}/{\sqrt{n}}>0$ and $\tau^*+{u_{p+2}}/{\sqrt{\lambda_n}}>0$. It is possible since $\tau^*>0$. Thus, combined with the assumption $\lambda_n\rightarrow+\infty$ as $n$ tends to $\infty$, the involved sequence tends to a strictly positive limit as $n$ tends to $\infty$. Since ${\lambda_n}/{n}\rightarrow0$ as $n$ tends to $\infty$, two cases are possible. Either, $\sqrt{{\lambda_n}/{n}}u_j\leq Lu_{p+2}$ and 
\begin{equation}
\label{L3}
b_{n,j}(u)=P_{n,j}(u)-\frac{\sqrt{\lambda_n}u_{p+2}} {\hat{w}^{adb}_j}
=\frac{\lambda_nu_j\hat{w}^{adb}_j}{\sqrt{n}},
\end{equation}
or $\sqrt{{\lambda_n}/{n}}u_j> Lu_{p+2}$ and
\begin{equation}
\label{L4}
b_{n,j}(u)=
\frac{\lambda_n\hat{w}^{adb}_j}{\left(\tau^*+\frac{u_{p+2}}{\sqrt{\lambda_n}}\right)}
\left(\frac{u_j^2}{2Ln}+\frac{\tau^*u_j}{\sqrt{n}}\right)+\frac{L\hat{w}^{adb}_j}{2\left(\tau^*+\frac{u_{p+2}}{\sqrt{\lambda_n}}\right)}u_{p+2}^2.
\end{equation}
Similiarly, we get the same result if  $\beta^*_j= -L\tau^*$. 
Gathering \eqref{nonL8} and using $B(\pm L)=0$, we have the following decomposition:
\begin{equation}
\label{eq8}
\sum_{j=1}^{p_0}P_{n,j}(u)=\sum_{j=1}^{p_0}c_{n,j}(u)+\sum_{j=1}^{p_0}\left(a_{n,j}(u)\indic{\vert\beta^*_j\vert\neq L\tau^*}+b_{n,j}(u)\indic{\vert\beta^*_j\vert= L\tau^*}\right)+\frac{{u_{p+2}}^2}{2L{\tau^*}^3}\sum_{j=1}^{p_0}\vert\beta^*_j\vert^{2-\gamma}\indic{\vert\beta^*_j\vert>L\tau*},
\end{equation}
where
$$c_{n,j}(u)=u_{p+2}\sqrt{\lambda_n}\sum_{j=1}^{p_0}\left(\frac 1 {\hat{w}^{adb}_j} +\hat{w}^{adb}_jB\left(\frac{\beta^*_j}{\tau^*}\right)\right).$$
We now study the convergence of each term.
The $\sqrt{n}$-consistency of $\hat{\beta}^{unpen}$ implies that $\hat{w}^{adb}_j \stackrel{\mathbb{P}}{\rightarrow}1/\vert\beta^*_j\vert^{\gamma}<+\infty$. Moreover, ${\lambda_n}/{\sqrt{n}}\rightarrow0$ as $n$ tends to infinity, thus, by Slutsky's theorem, the first three terms of $a_{n,j}(u)$ tends to $0$ in probability for any $(u)\in\rset^{p+2}$ fixed.  Concerning the last term (the rest), we have that 
$$\forall\,\eps>0\,,\exists\,N_{\eps}(u),\,\forall\,n\geq N_{\eps}(u),\,\lambda_n\xi\left(\frac{u_j}{\sqrt{n}},\frac{u_{p+2}}{\sqrt{\lambda_n}}\right)\leq\eps\left(\frac{u_j^2\lambda_n}{n}+u_{p+2}^2\right).$$ 
Moreover, $({\lambda_n}/{n})_{n\geq1}$ is a bounded sequence (since it converges to $0$ as $n$ tends to infinity). Thus, $\lambda_n\xi ({u_j}/{\sqrt{n}},u_{p+2}/{\sqrt{\lambda_n}})\rightarrow0$ as $n$ tends to $\infty.$ Consequently, for any $u\in\rset^{p+2}$ fixed, the forth term of $a_{n,j}$ tends to $0$ in probability. Using Slutsky's lemma, this entails that, for any $u\in\rset^{p+2}$ fixed,
$a_{n,j}(u)\stackrel{\mathbb{P}}{\rightarrow}0.$
Concerning the term $b_{n,j}(u)$
As previously we have $\hat{w}^{adb}_j \stackrel{\mathbb{P}}{\rightarrow}1/\vert\beta^*_j\vert^{\gamma}<+\infty$ and ${\lambda_n}/{\sqrt{n}}\rightarrow0$ as $n$ tends to infinity, so, if $\beta^*_j= L\tau^*$, 
$$
b_{n,j}(u)\stackrel{\mathbb{P}}{\rightarrow}\frac{L^{(1-\gamma)}}{2{\tau^*}^{(\gamma+1)}}u_{p+2}^2\indic{u_{p+2}<0}\,.
$$
Similarly, we get the same result if  $\beta^*_j= -L\tau^*$. Concerning the term $c_{n,j}(u)$,
Property \eqref{proptaustar} (see Lemma~\ref{existtaustar}) is available since $\beta^*\neq0$ and
\begin{equation}
\label{introd8}
c_{n,j}(u)=u_{p+2}\sqrt{\frac{\lambda_n}n}
\sum_{j=1}^{p_0}\left(\sqrt{n}(\vert\hat{\beta}^{unpen}_j\vert^{\gamma}-\vert\beta^*_j\vert^{\gamma})+B\left(\frac{\beta^*_j}{\tau^*}\right)\left(\sqrt{n}\frac{\left(\vert\hat{\beta}^{unpen}_j\vert^{\gamma}-\vert\beta^*_j\vert^{\gamma}\right)}{\vert\hat{\beta}^{unpen}_j\vert^{\gamma}\vert\beta^*_j\vert^{\gamma}}\right)\right).
\end{equation}
Since $\beta^*_j\neq0$, $x\rightarrow\vert x\vert^{\gamma}$ is differentiable at $\beta^*_j$ and  the Taylor-Young theorem entails that 
$$\sqrt{n}\left(\vert\hat{\beta}^{unpen}_j\vert^{\gamma}-\vert\beta^*_j\vert^{\gamma}\right)=\gamma\text{sign}(\beta^*_j)\vert\beta^*_j\vert^{\gamma-1}\sqrt{n}\left(\hat{\beta}^{unpen}_j-\beta^*_j\right)+\sqrt{n}\left(\hat{\beta}^{unpen}_j-\beta^*_j\right)\xi_j\left(\hat{\beta}^{unpen}_j\right)$$ with $\xi_j(x)\rightarrow0$ as $x$ tends to $\beta^*_j$. Now, the $\sqrt{n}$-consistency of $\hat{\beta}^{unpen}$ implies that the first term of this expansion is bounded in probability. It also entails that $\hat{\beta}^{unpen}_j\stackrel{\mathbb{P}}{\rightarrow}\beta^*_j$ which leads to $\xi_j\left(\hat{\beta}^{unpen}_j\right)\stackrel{\mathbb{P}}{\rightarrow}0$ since $\xi_j(x)\rightarrow0$ as $x$ tends to $\beta^*_j$. Consequently, the second term of this expansion is also bounded in probability and, finally,
$\sqrt{n}(\vert\hat{\beta}^{unpen}_j\vert^{\gamma}-\vert\beta^*_j\vert^{\gamma})=\CO_P(1).$
Since ${\lambda_n}/{n}\rightarrow0$ as $n$ tends to infinity, and $\vert\hat{\beta}^{unpen}_j\vert^{\gamma}\stackrel{\mathbb{P}}{\rightarrow}\vert\beta^*_j\vert^{\gamma}\neq0$, so $c_{n,j}(u)$ converges in probability to 0.
Combining \eqref{eq8} with all these convergences, the convergence in probability of \eqref{non08} is proved. 
Using first theorem 2.7 (vi) of \cite{VanderVaart1998} and then that convergence in probability is stronger than convergence in distribution (theorem 2.7 (ii) of \cite{VanderVaart1998}), we get that convergence in probability implies finite-dimensional convergence in \eqref{non08}. Theorem 5 of \cite{Knight97} implies that \eqref{non08} holds since the limit function $u\rightarrow u_{p+2}^2C(u_{p+2})$ is finite.

{\bf Step 2}. Next, we treat the sum of terms $P_{n,j}$ for $j>p_0,$ and first show that
\begin{equation}
\label{oui08}
\left(P_{n,p_0+1},\cdots,P_{n,p}\right)\rightarrow_{e-d}\left(I_{B_{p_0+1}},\cdots,I_{B_{p}}\right)\,,
\end{equation}
where $B_j=\{(u_{1:p},u_{p+2})\in\rset^{p+1},\,u_j=0\}$ and for a set $A$, $I_{A}$ denotes the indicator function of $A$ (i.e. $I_{A}(x)=0$ if $x\in A$ and $I_{A}(x)=+\infty$ otherwise). Let us put
\begin{equation}
\label{dd8}
q_{n,j}(u)=P_{n,j}(u)-\sqrt{\lambda_n}u_{p+2}\vert\hat{\beta}^{unpen}_j\vert^{\gamma}.
\end{equation}
Since $\hat{\beta}^{unpen}$ is a $\sqrt{n}$-consistent estimator and $j\in[p_0+1,p]$, 
$n^{\gamma/2}\vert\hat{\beta}^{unpen}_j\vert^{\gamma}$ is a tight sequence. 
Moreover, we have ${\lambda_n}/{n^{\gamma}}\rightarrow0$ as $n$ tends to infinity, thus $\forall u_{p+2}\in\rset$, 
$\sqrt{\lambda_n}u_{p+2}\vert\hat{\beta}^{unpen}_j\vert^{\gamma}=u_{p+2}\sqrt{{\lambda_n}{n^{-\gamma}}}(\sqrt{n}\vert\hat{\beta}^{unpen}_j\vert)^{\gamma}\cvp0.$
Using first theorem 2.7 (vi) of \cite{VanderVaart1998}, we get that convergence in probability implies finite-dimensional convergence:
$\sqrt{\lambda_n}u_{p+2}\vert\hat{\beta}^{unpen}_j\vert^{\gamma}\rightarrow_{f-d}0.$ Since the involved limit function is finite and by convexity,  theorem 5 of \cite{Knight97}  ensures  that we have the epiconvergence in distribution.
Moreover, $\IH{x}{L}\geq\vert x\vert$ and $\IH{0}{L}=0$, Lemma \ref{calcd2} with $q(x)=\IH{x}{L}$ leads to 
$$
d(q_{n,j},I_{B_j})\leq  2^{-\left[\tau^*\sqrt{\lambda_n}\right]+1}+\frac{2\sqrt{n}\vert\hat{\beta}^{unpen}_j\vert^{\gamma}}{\lambda_n},
$$
where $d$ is defined as in \eqref{defd}.
We have $\lambda_n\rightarrow+\infty$ as $n$ tends to infinity and $2^{-\left[\tau^*\sqrt{\lambda_n}\right]+1}\rightarrow0$ as $n$ tends to infinity since $\tau^*>0$. Furthemore ${2\sqrt{n}\vert\hat{\beta}^{unpen}_j\vert^{\gamma}}/{\lambda_n}={2(\sqrt{n}\vert\hat{\beta}^{unpen}_j\vert)^{\gamma}}/\lambda_n/n^{(\gamma-1)/2}$ and since $\hat{\beta}^{unpen}$ is a $\sqrt{n}$-consistent estimator and $j\in[p_0+1,p]$, the numerator is a tight sequence and the denominator tends to $+\infty$ as $n$ tends to infinity. Consequently, ${2\sqrt{n}\vert\hat{\beta}_j\vert^{\gamma}}/{\lambda_n}\cvp0$ and $d(q_{n,j},I_{B_j})\stackrel{\mathbb{P}}{\rightarrow}0$. 
Finally, using part (ii) of lemma 1.10.2 page 57 of \cite{VdWWellner96}, we have $q_{n,j}\rightarrow_{e-d}I_{B_j}.$
The notion of epi-convergence in distribution of convex lower semicontinuous random variables is a particular case of weak convergence of a net as stated in definition 1.33 of \cite{VdWWellner96}. Consequently, we can use Slutsky's theorem page 32, example 1.4.7 of \cite{VdWWellner96} to ensure that
\begin{equation}
\label{lims}
\left(\sqrt{\lambda_n}u_{p+2}\vert\hat{\beta}_j\vert^{\gamma},q_{n,j}(u_{1:p},u_{p+2})\right)\rightarrow_{e-d} (0,I_{B_j})
\end{equation}
since $0$ is deterministic. Moreover, we have $\sqrt{\lambda_n}u_{p+2}\vert\hat{\beta}^{unpen}_j\vert^{\gamma}\rightarrow_{u-d}0$ since 
we have shown the finite dimensional convergence in distribution and since $\sqrt{\lambda_n}u_{p+2}\vert\hat{\beta}_j\vert^{\gamma}$ and $0$ are finite convex functions (\cite{Arcones98} and \cite{Knight97}). We are now in position to use part (b) of theorem 4 of \cite{Knight97}: gathering \eqref{lims}, $\sqrt{\lambda_n}u_{p+2}\vert\hat{\beta}^{unpen}_j\vert^{\gamma}\rightarrow_{u-d}0$, continuity of $0$  and \eqref{dd8}, it ensures that $P_{n,j}\rightarrow_{e-d}I_{B_j}$ holds. Since $I_{B_j}$ is deterministic, theorem 18.10 (ii) of \cite{VanderVaart1998} ensures that the convergence in probability holds. Now, theorem 18.10 (vi) of  \cite{VanderVaart1998} leads to the convergence in probability in \eqref{oui08}. Moreover, convergence in probability is stronger than convergence in distribution thus \eqref{oui08} is proved.

For all $I\subset[p_0+1,p]$, $\text{dom}\left(\sum_{i\in I}I_{B_i}\right)=\{(u_{1:p},u_{p+2})\in \rset^{p+1},u_i=0,\forall\,i\in I\}$. Thus, for all $I\subset[p_0+1,p]$ and $J\subset[p_0+1,p]$ satisfying $I\cap J=\emptyset$,
$$
0\in \text{int}\left(\text{dom}\left(\sum_{i\in I}I_{B_i}\right)-\text{dom}\left(\sum_{j\in J}I_{B_j}\right)\right),
$$ 
where for $f$, a function defined on $\rset^{p+1}$, $\text{dom}(f)=\{x\in\rset^{p+1}/f(x)<+\infty\}$ and $A-B=\{a-b,\,a\in A,\,b\in B\}$.
Using successively this fact, \eqref{oui08}, Theorem 5 of \cite{McLindenBergstrom1981} and theorem 18.10 (iii) (v) (vi) and 18.11  of \cite{VanderVaart1998}, we get
\begin{equation}
\label{oui08s}
\sum_{j=p_0+1}^{p}P_{n,j}\rightarrow_{e-d}\sum_{j=p_0+1}^{p}I_{B_{j}}
\end{equation}
As previously, we can use Slutsky's theorem page 32, example 1.4.7 of \cite{VdWWellner96} to ensure that \eqref{oui08s}  and \eqref{non08} imply that 
\begin{equation}
\label{lims1}
\left(\sum_{j=p_0+1}^{p}P_{n,j}(u),\sum_{j=1}^{p_0}P_{n,j}(u)\right)\rightarrow_{e-d} \left(\sum_{j=p_0+1}^{p}I_{B_j},u_{p+2}^2C(u_{p+2})\right) 
\end{equation}
since $u_{p+2}^2C(u_{p+2})$ is deterministic. Moreover, we have $\sum_{j=1}^{p_0}P_{n,j}(u)\rightarrow_{u-d}u_{p+2}^2C(u_{p+2})$ since we have shown  the finite dimensional convergence in distribution and $\sum_{j=1}^{p_0}P_{n,j}(u)$ and $u_{p+2}^2C(u_{p+2})$ are finite (for $n$ sufficiently large) convex functions (\cite{Arcones98} and \cite{Knight97}). Using part (b) of theorem 4 of \cite{Knight97}: gathering 
\eqref{lims1},  $\sum_{j=1}^{p_0}P_{n,j}(u_{1:p},u_{p+2})\rightarrow_{u-d}{u_{p+2}}^2C(u_{p+2})$ and continuity of ${u_{p+2}}^2C(u_{p+2})$, it ensures that Lemma~\ref{cvpen4} holds. $
\blacksquare
$

\subsubsection{Proof of lemma~\ref{existtaustar}}

\begin{lemme}
\label{existtaustar}
If $\beta^*\neq0$ then there exists a unique $\tau^*>0$ satisfying equation (\ref{deftaustar}) and
\begin{equation}
\label{proptaustar}
\sum_{j=1}^{p}\vert\beta^*_j\vert^{\gamma}+\sum_{j=1}^{p_0}\frac{1}{\vert\beta^*_j\vert^{\gamma}}\left(\IH{\frac{\beta^*_j}{\tau^*}}{L}-\frac{\beta^*_j}{\tau^*}\dIH{\frac{\beta^*_j}{\tau^*}}{L}\right)=0.
\end{equation}
\end{lemme}
{\bf Proof.} Let us denote by $I$ the following function of $\tau$
$$I(\tau)=\tau \left( \sum_{j=1}^{p}\vert\beta^*_j\vert^{\gamma}+\sum_{j=1}^{p_0}\frac{1}{\vert\beta^*_j\vert^{\gamma}}\IH{\frac{\beta^*_j}{\tau}}{L}\right).$$
This function is convex and $I'(\cdot)$ is continuous, increasing with $I'(\tau)\rightarrow \sum_{j=1}^{p}\vert\beta^*_j\vert^{\gamma}$ as $\tau\rightarrow +\infty$ and, if $\beta^*\neq0$,  $I'(\tau)\rightarrow -\infty$ as $\tau\rightarrow 0$. This leads to the existence of $\tau*>0$ by the  intermediate value theorem. The minimum of $I$ is unique since $I'$ is strictly increasing on each pieces $]0,\vert\beta^*_{(1)}/L\vert[$ and 
$[\vert\beta^*_{(k}\vert/L,\vert\beta^*_{(k+1)}\vert/L[$ for $1\leq k\leq p-1$, continuous and increasing on $\rset_+^*$, strictly positive at $\vert\beta^*_{(p)}\vert/L$ since $I'(\vert\beta^*_{(p)}\vert/L)=\sum_{j=1}^{p}\vert\beta^*_j\vert^{\gamma}>0$. Note that $I'$ is constant on $[\vert\beta^*_{(p)}\vert/L,+\infty[$. This concludes the proof.    
$
\blacksquare
$

\subsubsection{Proof of lemma~\ref{calcd2}}

For $f$, a function defined on $S$, we note $epi(f),$ its epigraph given by $epi(f)=\{(x,t)\in S\times \rset / f(x)\leq t\}.$
\begin{lemme}
\label{calcd2}
Let $q$ be a function such that $q(0)=0$ and  $\forall\,x\in\rset$, $q(x)\geq\vert x\vert$ . We use the notations of the proof of lemma \ref{cvpen4}. Let us recall that $q_{n,j}(u_{1:p},u_{p+2})=P_{n,j}(u_{1:p},u_{p+2})-\sqrt{\lambda_n}u_{p+2}\vert\hat{\beta}_j\vert^{\gamma}$ where 
$$
P_{n,j}(u)=\left\{
\begin{array}{lll}
\lambda_n\left(\frac{u_{p+2}}{\sqrt{\lambda_n}}\vert\hat{\beta}_j\vert^{\gamma} + \frac{1}{\vert\hat{\beta}_j\vert^{\gamma}}\left(\frac{u_{p+2}}{\sqrt{\lambda_n}}+\tau^*\right)q\left(\frac {\frac{\gamma_j}{\sqrt{n}}+\beta^*_j}{\frac{u_{p+2}}{\sqrt{\lambda_n}}+\tau^*}\right)-\frac{\tau^*}{\vert\hat{\beta}_j\vert^{\gamma}}q\left(\frac{\beta^*_j}{\tau^*}\right)\right)&\text{if}& u_{p+2}>-\sqrt{\lambda_n}\tau^*,\\
-\lambda_n\tau^*\left(\vert\hat{\beta}_j\vert^{\gamma}+\frac{1}{\vert\hat{\beta}_j\vert^{\gamma}}q\left(\frac{\beta^*_j}{\tau^*}\right)\right)&\text{if}&  u_{p+2}=-\sqrt{\lambda_n}\tau^*,\\
&&\text{and}\,\gamma_j=-\sqrt{n}\beta^*_j,,\\
+\infty&&\text{otherwise}\,.\\
\end{array}
\right.
$$
Then, $\forall j\in [p_0+1,p]$, 
$$
d(q_{n,j},I_{B_j})\leq 2^{-\left[\tau^*\sqrt{\lambda_n}\right]+1}+\frac{2\sqrt{n}\vert\hat{\beta}_j\vert^{\gamma}}{\lambda_n},
$$
where
\begin{equation}
\label{defd}
d(q_{n,j},I_{B_j})=\sum_{k=1}^{+\infty}\frac{1\wedge d_k(\text{epi}(q_{n,j}),\text{epi}(I_{B_j}))}{2^k}\,,
\end{equation}
$d_k$ is a semi-distance ("constrained Pompeiu-Haussdorf distance") 
\begin{equation}
\label{defdk}d_k(\text{epi}(q_{n,j}),\text{epi}(I_{B_j}))=\underset{\|x\|\leq k}{\text{max}}\vert d_{\text{epi}(q_{n,j})}(x)-d_{\text{epi}(I_{B_j})}(x)\vert\,,\end{equation} 
and $d_S(x)=\underset{y\in S}{\text{min}}\|x-y\|$ for a subset $S$ of $\rset^{p+1}$.
\end{lemme}
{\bf Proof. }
Let us note that distance $d$ caracterises the epi-convergence of lower semi-continuous functions: a sequence $\{f_n\}$ of extended-real-valued lower semi-continuous functions from $\rset^{p+1}$ {\it epiconverges} to a extended-real-valued lower semi-continuous function $f$ if and only if $d(f_n,f)\rightarrow 0$ as $n$ goes ton infinity.
We recall that $B_j=\{(u_{1:p},u_{p+2})\in\rset^{p+1},\,u_j=0\}$ and for a set $A$, $I_{A}$ denotes the indicator function of $A$. Let us introduce the set $D_j=\{(u_{1:p},u_{p+2})\in\rset^{p+1},\,u_j=0\,\text{and}\,u_{p+2}\geq-\sqrt{\lambda_n}\tau^*\}$. By using the triangular inequality,
\begin{equation}
\label{dec001}
d(q_{n,j},I_{B_j})\leq d(q_{n,j},I_{D_j})+d(I_{D_j},I_{B_j})\,.
\end{equation}
To begin with, let us show that 
\begin{equation}
\label{term001}
d(I_{D_j},I_{B_j})\leq 2^{-\left[\tau^*\sqrt{\lambda_n}\right]}\,.
\end{equation}
Here we use a geometrical point of view. The epigraph of the indicator function $I_{A}$ of a set $A$ is the ``half- cylinder with cross-section $A$'' i.e. $A\times\rset_+$. Consequently, the  epigraph of $I_{B_j}$ is an half-hyperplan supported by the $u_j$ axis and the  epigraph of $I_{D_j}$ is the part of this half-hyperplan where, moreover, $u_{p+2}\geq-\sqrt{\lambda_n}\tau^*$. Note that this cut is perpendicular to the $u_{p+2}$-axis. So if we consider $x\in\rset^{p+2}$ such that $x_{p+1}\geq-\sqrt{\lambda_n}\tau^*$, the distance between $x$ and $\text{epi}(I_{D_j})$ is reached for a point in $\text{epi}(I_{B_j})$. Thus
\begin{equation}
\label{epieq1}
\forall\,x,\,\|x\|_2\leq k\, \text{with}\,k\leq \sqrt{\lambda_n}\tau^*,\, d_{\text{epi}(I_{D_j})}(x)=d_{\text{epi}(I_{B_j})}(x), 
\end{equation}
and if $k\leq\sqrt{\lambda_n}\tau^*$ then $d_k(\text{epi}(I_{D_j}),\text{epi}(I_{B_j}))=0.$ Now the definition \eqref{defd} of the distance $d$ implies that
$$
d(I_{D_j},I_{B_j})=\sum_{k\geq\left[\sqrt{\lambda_n}\tau^*\right]+1}\frac{1\wedge d_k(\text{epi}(I_{D_j}),\text{epi}(I_{B_j}))}{2^k}\leq \sum_{k\geq\left[\sqrt{\lambda_n}\tau^*\right]+1}\frac{1}{2^k}\,,
$$
and \eqref{term001} is  proved.

Next, we show that 
\begin{equation}
\label{term002}
d(q_{n,j},I_{D_j})\leq\frac{2\sqrt{n}\vert\hat{\beta}^{unpen}_j\vert^{\gamma}}{\lambda_n}+2^{-\left[\tau^*\sqrt{\lambda_n}\right]}\,.
\end{equation}
For $j\in[p_0+1,p]$, $q(0)=0$ implies that
\begin{equation}
\label{exprQ}
q_{n,j}(u_{1:p},u_{p+2})=\frac{\lambda_n}{\vert\hat{\beta}^{unpen}_j\vert^{\gamma}}\left(\frac{u_{p+2}}{\sqrt{\lambda_n}}+\tau^*\right)q\left(\frac{u_j}{\sqrt{n}\left(\frac{u_{p+2}}{\sqrt{\lambda_n}}+\tau^*\right)}\right)+I_{E}
\end{equation}
where we set $0/0=0$ and  $$E=\{(u_{1:p},u_{p+2}),\,u_{p+2}>-\sqrt{\lambda_n}\tau^*\}\cup\{(u_{1:p},u_{p+2}),\,u_{p+2}=-\sqrt{\lambda_n}\tau^*\,\text{and}\,u_j=0\}.$$ Consequently, $q_{n,j}(u_{1:p},u_{p+2})\leq I_{D_j}(u_{1:p},u_{p+2})$. Indeed, it is clear if $(u_{1:p},u_{p+2})\notin D_j$. Moreover, if $(u_{1:p},u_{p+2})\in D_j$, $q_{n,j}(u_{1:p},u_{p+2})=0$ since $q(0)=0$. Consequently, $\text{epi}\left(I_{D_j}\right)\subset\text{epi}\left(q_{n,j}\right)$, $d_{\text{epi}\left(I_{D_j}\right)}(.)\geq d_{\text{epi}\left(q_{n,j}\right)}(.)$ and 
\begin{equation*}
d_k(\text{epi}\left(q_{n,j}\right),\text{epi}\left(I_{D_j}\right))=\underset{\|x\|\leq k}{\text{max}}\left( d_{\text{epi}(I_{D_j})}(x)-d_{\text{epi}(q_{n,j})}(x)\right)\,.
\end{equation*}
Since  $\forall\,t\in\rset$, $q(t)\geq\vert t\vert$, it holds that,$\forall\,(t,\tau)\in\rset\times\rset_+^*$,  $\tau q({t}/{\tau})\geq\vert t\vert$ and expression \eqref{exprQ} entails 
$$
q_{n,j}(u_{1:p},u_{p+2})\geq F_{n,j}(u_{1:p},u_{p+2})\,,
$$
where $F_{n,j}(u_{1:p},u_{p+2})=
{\lambda_n\vert u_j\vert\vert\hat{\beta}^{unpen}_j\vert^{-\gamma}}/{\sqrt{n}}+I_{E}$. Consequently, $\text{epi}\left(q_{n,j}\right)\subset\text{epi}\left(F_{n,j}\right)$, $d_{\text{epi}\left(q_{n,j}\right)}(.)\geq d_{\text{epi}\left(F_{n,j}\right)}(.)$ and 
\begin{equation}
\label{dk1}
d_k(\text{epi}\left(q_{n,j}\right),\text{epi}\left(I_{D_j}\right))\leq\underset{\|x\|\leq k}{\text{max}}\left( d_{\text{epi}(I_{D_j})}(x)-d_{\text{epi}(F_{n,j})}(x)\right)\,.
\end{equation}
Now, $
\text{epi}(F_{n,j})=S_1\cup S_2$ 
where $$S_1=\{(u_{1:p},u_{p+2},t)\in\rset^{p+2},\,u_{p+2}>-\sqrt{\lambda_n}\tau^*\,\text{and}\,\frac{\lambda_n\vert u_j\vert}{\sqrt{n}\vert\hat{\beta}^{unpen}_j\vert^{\gamma}}\leq t\},$$ and $$S_2=\{(u_{1:p},u_{p+2},t)\in\rset^{p+2},\,u_{p+2}=-\sqrt{\lambda_n}\tau^*\,,u_j=0,\text{and}\,t\geq0\}.$$ Thus,
\begin{equation}
\label{exprdepiF}
d_{\text{epi}(F_{n,j})}(x)=d_{S_1}(x)\wedge d_{S_2}(x)\,.
\end{equation}
Easy calculations lead to, $\forall\,x\in\rset^{p+2}$,
\begin{equation}
\label{exprd2}
d^2_{S_2}(x)=\inf_{z\in S_2}\sum_{i=1}^{p+2}(x_i-z_i)^2=x_j^2+(x_{p+1}+\sqrt{\lambda_n}\tau^*)^2+x_{p+2}^2\indic{x_{p+2}<0},
\end{equation}
and
$$
d^2_{S_1}(x)=\inf_{z\in S_2}\sum_{i=1}^{p+2}(x_i-z_i)^2=d^2_{\text{epi}(f_{n,j})}(x_1,\cdots,x_p,x_{p+2})+(x_{p+1}+\sqrt{\lambda_n}\tau^*)^2\indic{x_{p+1}<-\sqrt{\lambda_n}\tau^*},
$$
where $f_{n,j}(u_{1:p})={\lambda_n}\vert u_j\vert\vert\hat{\beta}^{unpen}_j\vert^{-\gamma}/\sqrt{n}$.
If we consider $x\in\rset^{p+2}$ such that $\|x\|_2\leq k$ with $k\leq\sqrt{\lambda_n}\tau^*$, it satisfies that $x_{p+1}\geq-\sqrt{\lambda_n}\tau^*$ and thus $d^2_{S_1}(x)=d^2_{\text{epi}(f_{n,j})}(x_1,\cdots,x_p,x_{p+2})\,.$ Technical computations leads to
\begin{equation}
\label{df2}
d_{S_1}(x)=\left\{
\begin{array}{lll}
\sqrt{x_j^2+x_{p+2}^2}&\text{if}&x_{p+2}\leq-\frac{\sqrt{n}\vert\hat{\beta}^{unpen}_j\vert^{\gamma}}{\lambda_n}\vert x_j\vert,\\
\frac{\vert x_j\vert-x_{p+2}\frac{\sqrt{n}\vert\hat{\beta}^{unpen}_j\vert^{\gamma}}{\lambda_n}}{\sqrt{1+\frac{n\vert\hat{\beta}^{unpen}_j\vert^{2\gamma}}{\lambda_n^2}}}&\text{if}& -\frac{\sqrt{n}\vert\hat{\beta}^{unpen}_j\vert^{\gamma}}{\lambda_n}\vert x_j\vert<x_{p+2}\leq \frac{\lambda_n}{\sqrt{n}\vert\hat{\beta}^{unpen}_j\vert^{\gamma}}\vert x_j\vert,\\
0&\text{if}& x_{p+2}> \frac{\lambda_n}{\sqrt{n}\vert\hat{\beta}^{unpen}_j\vert^{\gamma}}\vert x_j\vert.
\end{array}
\right.
\end{equation}
Using explicit expressions \eqref{df2} and \eqref{exprd2}, we can show that for any $x\in\rset^{p+2}$ such that $\|x\|_2\leq k$ with $k\leq\sqrt{\lambda_n}\tau^*$, 
\begin{equation}
\label{in000}
d_{S_1}(x)\leq d_{S_2}(x).
\end{equation}
Gathering \eqref{in000} with \eqref{exprdepiF}, for any $x\in\rset^{p+2}$ such that $\|x\|_2\leq k$ with $k\leq\sqrt{\lambda_n}\tau^*$, 
\begin{equation}
\label{epieq2}
d_{\text{epi}(F_{n,j})}(x)=d_{S_1}(x)=d_{\text{epi}(f_{n,j})}(x_1,\cdots,x_p,x_{p+2}).
\end{equation}
Combining \eqref{dk1}, \eqref{epieq2} and \eqref{epieq1}, if $k\leq \sqrt{\lambda_n}\tau^*$, we obtain
$$
d_k(\text{epi}\left(q_{n,j}\right),\text{epi}\left(I_{D_j}\right))\leq\underset{\|x\|\leq k}{\text{max}}\left( d_{\text{epi}(I_{B_j})}(x_1,\cdots,x_{p},x_{p+2})-d_{\text{epi}(f_{n,j})}(x_1,\cdots,x_{p},x_{p+2})\right).
$$
The involved objective function does not depend on $x_{p+1}$. Moreover, using the form of the constraints, if $k\leq \sqrt{\lambda_n}\tau^*$, we get 
$$
d_k(\text{epi}\left(q_{n,j}\right),\text{epi}\left(I_{D_j}\right))\leq\underset{x_1^2+\cdots+x_{p}^2+x_{p+2}^2\leq k^2}{\text{max}}\left( d_{\text{epi}(I_{A_j})}(x_1,\cdots,x_{p},x_{p+2})-d_{\text{epi}(f_{n,j})}(x_1,\cdots,x_{p},x_{p+2})\right).
$$
Moreover, since $\forall u_{1:p}\in\rset^{p}$, $I_{A_j}(u_{1:p})\geq f_{n,j}(u_{1:p})$, if $k\leq \sqrt{\lambda_n}\tau^*$,
$$
d_k(\text{epi}\left(q_{n,j}\right),\text{epi}\left(I_{D_j}\right))\leq d_k(\text{epi}\left(f_{n,j}\right),\text{epi}\left(I_{A_j}\right)),
$$
and technical computations leads to
$$
d_k(\text{epi}\left(f_{n,j}\right),\text{epi}\left(I_{A_j}\right))=\frac{k\sqrt{n}\vert\hat{\beta}^{unpen}_j\vert^{\gamma}}{\lambda_n\sqrt{1+\frac{n\vert\hat{\beta}^{unpen}_j\vert^{2\gamma}}{\lambda_n^2}}}\,.
$$
Finally, using the definition \eqref{defd}, we have
$$
d\left(q_{n,j},I_{D_j}\right)\leq \sum_{k\leq \left[\sqrt{\lambda_n}\tau^*\right]}\frac{d_k(\text{epi}\left(q_{n,j}\right),\text{epi}\left(I_{D_j}\right))}{2^k}+\sum_{k\geq \left[\sqrt{\lambda_n}\tau^*\right]+1}\frac{1}{2^k}\,.
$$ 
Gathering this inequality with the previous one and the fact that $\sum_{k\geq1}\frac{k}{2^k}\leq2$, \eqref{term002} is proved. Using
equation \eqref{dec001} with  \eqref{term001} and \eqref{term002}, the bound involved in Lemma~\ref{calcd2} holds.
$\blacksquare$

\section*{Acknowledgements}
Part of this work was supported by the Interuniversity Attraction
Pole (IAP) research network in Statistics P5/24 and by MSTIC project of the Joseph-Fourier University. We are grateful to
Anestis Antoniadis for constructive and fruitful discussions. 


\newpage
\section*{ Tables and Figures}

\begin{table}[h]
\caption{Selection model ability on Model 1 based on 100 replications.}
\label{modselecMod1}
\begin{center}
\begin{tabular}{lrrrrrrrr}
&C&O&U&Z&CZ&TZ&CNZ&TNZ\\ \hline 
\multicolumn{4}{l}{\bf Least square criterion, $n=100$}\\
{\tt ad-lasso} & 0 & 0 & 100 & 35.38& 24.32 & 25&3.94 &15\\
{\tt ad-en}  & 0 & 0 & 100 & 27.72  &19.69 & 25&6.97&15\\
{\tt ad-Berhu} & 0 & 51 & 49 & 6.50  &5.79 &25& 14.27&15\\
\multicolumn{4}{l}{\bf Huber's criterion, $n=100$}\\
{\tt ad-lasso} & 0 & 3 & 97 &20.9  &14.08 &25&8.18&15\\
{\tt ad-en} & 0 & 3 & 97 & 26.99  & 18.81 &25& 6.82&15\\
{\tt ad-Berhu} & 0 & 33 & 67 & 9.89  & 8.57 &25& 13.68 &15\\
\hline 
\multicolumn{4}{l}{\bf Least square criterion, $n=200$}\\
{\tt ad-lasso} &0  &0  & 100& 35.15  & 24.59 &25&4.44&15\\
{\tt ad-en} & 0 & 2 & 98 &27.71  &20.52 &25&7.81&15\\
{\tt ad-Berhu} & 0 & 42 & 58 &18.36 &17.31&25&13.95&15\\
\multicolumn{4}{l}{\bf Huber's criterion, $n=200$}\\
{\tt ad-lasso} & 0 & 7 &  93 & 29.52  & 19.85  & 25&5.33&15\\
{\tt ad-en} & 0 & 1 & 99 & 30.41  & 21.08 &25& 5.67&15\\
{\tt ad-Berhu} &0  & 48 & 52 &16.77  & 15.86&25&14.09&15\\ 
\hline 
\multicolumn{4}{l}{\bf Least square criterion, $n=400$}\\
{\tt ad-lasso} & 0 & 0 & 100 & 34.75 &24.74 &25&4.99&15\\
{\tt ad-en} & 0 & 0 & 100 & 27.06  & 21.17 &25& 9.11&15\\
{\tt ad-Berhu} & 1 & 32 & 67& 21.46 &19.29 &25&12.83&15\\ 
\multicolumn{4}{l}{\bf Huber's criterion, $n=400$}\\
{\tt ad-lasso}  & 0 & 9 & 91  & 29.67  & 20.3 &25& 5.63&15 \\
{\tt ad-en}  & 0 & 0 & 100 & 30.01 & 21.59 &25& 6.58&15\\
{\tt ad-Berhu} & 0 & 35 & 65  & 19.47  & 18.42 &25& 13.95&15\\ 
\hline 
\end{tabular}
\end{center}
\end{table}

\begin{table}[h]
\caption{Selection model ability on model 2 based on 100 replications.}
\label{modselecMod2}
\begin{center}
\begin{tabular}{lrrrrrrrr}
&C&O&U&Z&CZ&TZ&CNZ&TNZ\\ \hline 
\multicolumn{4}{l}{\bf Least square criterion, $n=100$}\\
{\tt ad-lasso} & 0 & 0 & 100&35.03 & 24.18&25&4.15&15\\ 
{\tt ad-en} & 0 & 0 & 100 & 28.66 & 20.46 &25& 6.80&15\\ 
{\tt ad-Berhu}& 0 & 47 & 53& 8.02 & 7.14&25&14.12&15\\ 
\multicolumn{4}{l}{\bf Huber's criterion, $n=100$}\\
{\tt ad-lasso} & 0 & 22 & 78 & 10.79  & 7.38&25  & 11.59&15\\ 
{\tt ad-en} & 1 & 15 & 84 &20.96  &16.09 &25&10.13 &15\\ 
{\tt ad-Berhu}& 5 & 24 & 71 & 24.36  & 23.12 &25& 13.76&15\\ 
\hline 
\multicolumn{4}{l}{\bf Least square criterion, $n=200$}\\
{\tt ad-lasso} & 0 & 0 & 100&35.25  &24.61 &25&4.36&15\\ 
{\tt ad-en} & 0 & 0  & 100 & 28.21  & 20.68 &25& 7.47&15\\ 
{\tt ad-Berhu}& 0 & 49 & 51 &17.17 &16.59 &25& 14.11&15\\ 
\multicolumn{4}{l}{\bf Huber's criterion, $n=200$}\\
{\tt ad-lasso} & 0 & 6  & 94 & 21.63 & 16.77 &25& 10.14&15\\ 
{\tt ad-en} & 1 & 15 & 84 &25.38  &20.31 &25&9.93 &15\\ 
{\tt ad-Berhu}& 4 & 11 & 85 & 25.36 &23.78 &25& 13.42&15\\ 
\hline 
\multicolumn{4}{l}{\bf Least square criterion, $n=400$}\\
{\tt ad-lasso} & 0 & 0 & 100 &35.05  &24.69 &25&4.64&15\\ 
{\tt ad-en} & 0 & 0 & 100 & 28.10   &21.69 &25&8.59&15\\ 
{\tt ad-Berhu} & 1 & 30 & 69 &20.15  &18.74 &25&13.59&15\\ 
\multicolumn{4}{l}{\bf Huber's criterion, $n=400$}\\
{\tt ad-lasso} & 0 & 14 & 86 & 21.30 &17.76 &25&11.46&15\\ 
{\tt ad-en} & 2 & 25 & 73 &23.32  &19.90 &25&11.58&15\\ 
{\tt ad-Berhu} & 23 & 10 & 67 &24.64  & 23.59  &25& 13.95&15\\ 
\hline 
\end{tabular}
\end{center}
\end{table}

\begin{table}[h]
\caption{Selection model ability on model 3 based on 100 replications.}
\label{modselecMod3}
\begin{center}
\begin{tabular}{lrrrrrrrr}
&C&O&U&Z &CZ&TZ&CNZ&TNZ\\ \hline 
\multicolumn{4}{l}{\bf Least square criterion, $n=100$}\\
{\tt ad-lasso}  & 0 & 0 & 100& 35.23 &24.45 &25& 4.22&15\\ 
{\tt ad-en}  &0  &0  &100 &29.11  &20.97 &25&6.86&15\\ 
{\tt ad-Berhu}  & 0 & 53 & 47 &6.68  &6.02&25 &14.34&15\\ 
\multicolumn{4}{l}{\bf Huber's criterion, $n=100$}\\
{\tt ad-lasso} & 0 & 4 & 96 &25.5  &17.36 &25&6.86 &15\\  
{\tt ad-en}  & 0 & 0 &100 &29.40  &20.60 &25&6.20 &15\\ 
{\tt ad-Berhu}  & 0 & 34 & 66 &14.34  &13.24 &25 &13.90 &15\\ 
\hline 
\multicolumn{4}{l}{\bf Least square criterion, $n=200$}\\
{\tt ad-lasso}  & 0 & 0 &100 &35.32  &24.74&25 &4.42&15\\ 
{\tt ad-en}  & 0 & 0  & 100 &28.58  &20.83 &25&7.25&15\\ 
{\tt ad-Berhu}  & 0 & 44 & 56 & 19.01 & 17.47&25&13.46&15\\ 
\multicolumn{4}{l}{\bf Huber's criterion, $n=200$}\\
{\tt ad-lasso}  & 0 & 3 & 97 &31.95  & 21.85 &25&4.9&15\\ 
{\tt ad-en}  & 0 & 2 &98 &28.37  &19.91 &25&6.54&15\\ 
{\tt ad-Berhu}  & 0 & 50  & 50  &19.53 &18.73&25 & 14.20&15\\ 
\hline 
\multicolumn{4}{l}{\bf Least square criterion, $n=400$}\\
{\tt ad-lasso}  & 0 & 0 &100 &34.78 &24.74&25 &4.96&15\\ 
{\tt ad-en}  & 0 & 0 &100 &26.54  &20.61 &25&9.07&15\\ 
{\tt ad-Berhu}  & 0 & 22 & 78&20.88  &19.29&25 &13.41&15\\ 
\multicolumn{4}{l}{\bf Huber's criterion, $n=400$}\\
{\tt ad-lasso}  & 0 & 7 &93 &31.68   &22.04 &25&5.36 &15\\ 
{\tt ad-en}  & 0 & 3 & 97 &26.25  &19.65 &25&8.4 &15\\ 
{\tt ad-Berhu}  & 2 & 23 & 75 &22.14   &21.02 &25&13.88 &15\\ 
\hline 
\end{tabular}
\end{center}
\end{table}

\begin{table}[h]
\caption{Prostate cancer data: comparing methods}
\label{prederr}
\begin{center}
\begin{tabular}{lll}
Methods&mean of 100 parameters (std of the 100) &mean of 100 RPE (std of the 100 )\\ \hline 
OLS  & none  & 0.6054(0.1397) \\
\multicolumn{3}{c}{\bf Least square criterion}\\
{\tt ad-lasso} & $\lambda_{n}:$ 2.4177(1.7368) & 0.6357(0.1410) \\
{\tt ridge}  &  $\lambda_{n}:$ 2.6104(2.3111) & 0.6145(0.1406)\\
{\tt ad-en}  & $\lambda_{1,n}:$ 1.1361(1.0048), $\lambda_{2,n}:$ 2.5032(10.2605)  & 0.6231(0.1351)\\
{\tt ad-Berhu}  &$\lambda_{n}:$ 1.9850(1.2782)  &0.6237(0.1423) \\
\multicolumn{3}{c}{\bf Huber's criterion}\\
{\tt ad-lasso}  & $\lambda_{n}:$ 26.2749(7.4369) &0.7765(0.1879) \\
{\tt ridge}  & $\lambda_{n}:$ 3.7437(3.5792) & 0.6020(0.1327) \\
{\tt ad-en}  &  $\lambda_{1,n}:$ 1.3885(1.5778), $\lambda_{2,n}:$ 4.3222(14.2073) & 0.6185(0.1295)\\
{\tt ad-Berhu}  &$\lambda_{n}:$ 2.7456(1.9015)  & 0.6322(0.1391)\\
\end{tabular}
\end{center}
\end{table}

\end{document}